\documentclass{article}
\usepackage{amsmath}
\usepackage[english]{babel}
\usepackage[utf8]{inputenc}
\usepackage{graphicx}
\usepackage{soul}
\usepackage{epstopdf}
\usepackage[bottom=2.5cm,left=3cm,right=3cm,top=3cm]{geometry}

\begin{document}


\begin{center}
\section*{\LARGE Discrete two-countour system with one-directional motion}

\textit{Yashina M.V.$^{1,2}$, Tatashev A.G.$^{1,2}$, Fomina M.J.$^{1}$}\\
\textit{$^{1}~$Moscow Automobile and Road Construction State Technical University (MADI)}\\
\textit{$^{2}~$Moscow Technical University of Communications and Informatics (MTUCI)}\\
e-mail:mv.yashina@madi.ru
\end{center}

\large
{\bf Abstract.} This paper studies a dynamical system, which contains
two contours. There is a cluster on each contour. The cluster
contains particles, located in adjacent cells. The clusters
move under prescribed rules.
The delays of clusters are due to that the clusters
cannot pass through the node simultaneously.
The dynamical system belongs to the class of contour networks introduced
by A.P.~Buslaev.

\section*{1. Introduction} In versions of Nagel--Schreckenberg model, [1],
the particles move on a closed or infinite one-dimensional lattice under prescribed rule. In of these versions, at any discrete moment, each particle
moves onto one cell in the direction of movement if the cell ahead
is vacant. This rule is elementary cellular automaton CA~184 in the classification of S.~Wolfram, [2]. In [3[--[7], analytical results are obtained for one-contour models of this type.

A two-dimensional traffic model with a toroidal supporter (BML traffic model)
has been introduced in [7]. Conditions
of self-organization (all particles move without delays after some moment) and collapse (no particle moves after some moment)
for BML model have been obtained in [8]--[10].

In [12], the concept of a
traffic model with cluster movement has been introduced. In the discrete
version of the cluster model, each contour contains a given
quantity of cells. There are clusters of particles on each contour.
Particles of each cluster occupy neighboring cells. All particles of each
cluster move simultaneously in accordance with the rule of the cellular
automaton~240. Clusters can be delayed at nodes. In the continuous version of the model a cluster is a segment moving on the contour with constant velocity in a given direction.

The concept of a contour network has been introduced by A.P. Buslaev, [13].
In papers A.P.~Buslaev, basic approaches
for the study of complex networks has been developed. In accordance
with these approaches, models are dynamical systems.
The supporters of these systems are contours systems
with network structures. Particles (clusters) move
on contours in accordance with some rules. There are common nodes of contours called nodes. Delays in movement at nodes are due to that two contours cannot move through the same node simultaneously.
Analytical results have been obtained for systems with
two contours and for systems with regular periodic
structures and any number of contours.

Analytical results for contour networks were obtained in [13]--[23].
In partucular, in [22], [23] a two-contour system with two nodes is studied under the assumtion that the length of both clusters is the same. Two versions of the system were considered. In the system with co-directional motion one cluster moves counter-clockwise and the  other cluster moves clockwise. In the system with one-directional motion  both the clusters move counter-clockwise.

In this paper we study a discrete two-contour system with two nodes and clusters of different lengths. We consider the case of one-dimensional motion.

\section*{2. System description}

\hskip 18pt Suppose that the system contains two closed contours~----{\it contour~1}  and {\it contour~2}, Fig. 1. On contour $i$ there are the same number $n$ of {\it cells} and {\it cluster}  containing $l_i <n$ {\it particles}, which are located in $l_i$ neighboring cells and move simultaneously at discrete times $t = 0, 1,2, \dots $, $i = 1, 2.$ The cells have numbers $0, 1,\dots, n - 1$ increasing in the direction of cluster motion module $n.$ There are two common points of the contours, called {\it nodes}. The {\it node 1} is located between the cells $0$ and $1$ on the contour 1 and between the cells $d$, $d+1$ $(d\le n/2)$ on the contour 2. The {\it node 2} is located between cells $0$ and $1$ on contour 2 and between cells $d$, $d+1$ on contour 1.
If there is no delay, all particles of the cluster move one cell in the direction of motion at each discrete moment. The delays are caused by the impossibility of simultaneous passage of particles of different contours through the node.  We suppose the cluster of contour $i$ {\it(cluster $i$) is at node $i$} at time $t$, if at this moment the front particle of this cluster is in cell 0, and cluster $i$ is {\it at node} $j\ne i$, if the front particle of the cluster is in cell $d$, $i = 1, 2$. We suppose a cluster of contour $i$ {\it (cluster $i$) occupies node $i$} at time $t$, if at this moment there are particles in cells 0, 1 of this contour and {\it occupies node $j \ne i)$} if at the moment there are particles in cells $d$, $d+1$ of contour $i$, $i = 1, 2$. If at time $t$ cluster $i$ is at one of the nodes, and another cluster occupies this node, then there is a delay in the movement of cluster $i$ and in the moment $t+1$ the particles of this cluster are in the same cells as at the moment $t$. If the clusters are located at the same time at the node $i$ ({\it a competition occurs),} then the cluster $j\ne i$
 passes through the node first. In accordance with this {\it competition resolution rule,} the cluster passing to the arc of not less length {\it wins the competition.}

\begin{figure}[ht!]
\centerline{\includegraphics[width=300pt]{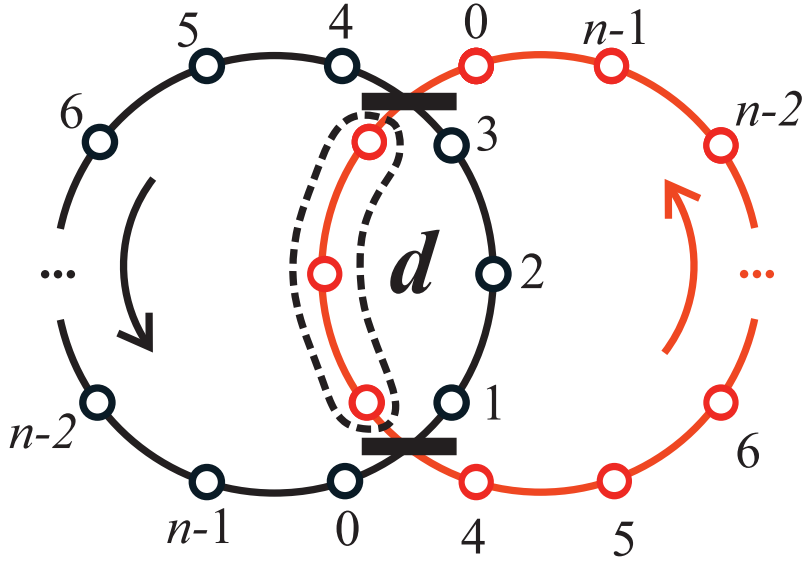}}
\caption{A two-contour system}
\end{figure}

We suppose that $d\le \frac{1}{2},$ $l_1\le l_2,$ which, due to symmetry, does not limit the generality.

The state of the system at time $t$ is a vector $A(t)=(\alpha_1(t),\alpha_2(t)),$ where $\alpha_i(t)$~--- the number of the cell in which the front particle of the cluster $i$, $ i = 1,2$ is located.  The state of the system is {\it acceptable,} if, when the system is in this state, the clusters do not occupy the same node at the same time. {\it Initial state,} i.e. the state at the initial moment of time $ t = 0, $ is specified. This state must be valid. It is clear that with the given rules of motion, the system can only go from valid state to valid one.

\section*{3. Limit cycles. Average cluster speed. States of free motion and collapse}

 \hskip 18pt Since the number of states is finite and the behavior of the system is deterministic, starting from a certain moment its states are repeated after $ T $ cycles. We have a {\it limit cycle} with period $ T. $

Let $A_i$~--- the number of transitions of the cluster $ i $ on the limit cycle with the period $ T. $ Then the value $ v_i = A_i / T $ is called the average speed of the cluster $ i, $ $ i = 1,2. $

The average cluster velocities generally depend on the values of $ n, $ $ l_1, $ $ l_2, $ $ d_1, $ $ d_2 $ and the initial state of the system.

If at any moment $ t \ge t_0 $ all particles move, then we say that from the moment $ t_0 $ the system {\it is in a state of free motion.}
We will call {\it self-organization} the property of a system to enter a state of free movement from any initial state. If the system is in a state of free motion, then $ v_1 = v_2 = 1. $

If starting from the moment $ t_0 $ no particle moves, then we say that from the moment $ t_0 $ the system is in the state of {\it collapse.} If the system is in the state of collapse, then $ v_1 = v_2 = 0. $

 \section*{3. Some lemmas}

\hskip 18ptLet us prove some lemmas.
\vskip 3pt
{\bf Lemma 1.} {\it If the inequality $ l_1 + l_2> n$ is satisfied the system cannot be in a state of free motion.
\vskip 3pt
 Proof.} If at the limit cycle the system is in a state of free motion, then for a cycle period equal to $ n, $ each node is crossed by $ l_1 $ of particles of cluster~1 and $ l_2 $ of particles of cluster~2. But this is impossible when the inequality $ l_1 + l_2> n$ is satisfied,  since the same node cannot simultaneously cross more than one particle. Lemma~1 is proved.
\vskip 3pt
{\bf Lemma 2.} {\it For any limit cycle, the average velocities of both clusters are not equal to 0, or the system is in a state of collapse.
\vskip 3pt
Proof.} If one cluster does not move on a spectral cycle, then it is at one of the nodes, which is occupied by another cluster throughout the cycle. But this is impossible, since, under the assumptions of the lemma, the second cluster moves and must completely traverse the contour at least once during the limit cycle. This contradiction proves the lemma.

\vskip 3pt
{\bf Lemma 3.} {\it If at the limit cycle the system is not in a state of free motion or in a state of collapse, then the system falls on this limit cycle in one of the states $(l_1+d,0),$ $(0,l_2+d),$ $(l_1+d-n,0),$ $(0,l_2+d-n),$ $(l_1,d),$  $(d,l_2)$
\vskip 3pt
 Proof.}It follows from Lemma 2 and terms of Lemma 3 that on the considered limit cycle, the average velocities of both clusters are not equal, and at least one cluster is delayed. If there is a cluster~1 delay at node~1, then at the end of this delay the system is in the state $ (0, l_2 + d) $ (addition module $ n$).  If there is a cluster~1 delay at node~2, then this delay ends when the system is in state $ (d, l_2). $ If there is a delay of cluster~2 at node~1, then at the end of this delay the system is in state $ (l_1 + d, 0) $ (addition over module $ n). $ If there is a delay of the cluster~2 at the node~2, then this delay ends when the system is in the state $ (l_l, d). $ It follows that the system falls on the limit cycle in at least one of terms of the lemma~3 state. Lemma~3 is proved.
\vskip 3pt
{\bf Lemma 4.} {\it collapse state exists if and only if terms
$l_1>d$, $l_2>d$ are  met.
\vskip 3pt
Proof.} If terms of Lemma 4 are met, the collapse state is $(d,d).$

If the system is in a state of collapse, then each cluster is at one node and occupies another, but this is impossible under the terms of Lemma~4, since if this terms is not met, at least one of the clusters cannot simultaneously be at one node and occupy another node. Lemma~4 is proved.

 \section*{4.System behavior options}

 \subsection*{4.1. Self-organization on terms $l_2\le d,$ $l_1+l_2\le n-2d$ }

 \hskip 18pt We prove a theorem that gives a sufficient terms for self-organization of the system.

 \vskip 3pt
 {\bf Theorem 1.} {\it If terms
 $$l_2 \le d,$$
 $$l_1+l_2\le n-2d,$$
 are met, then the system is in a state of free movement at any initial state from a certain point in time.
 \vskip 3pt
 Proof.} Under the assumptions of Theorem 1, the system, in accordance with Lemma~4, cannot be in a state of collapse. If on some limit cycle the system is also not in a state of free motion, then in accordance with Lemma~3 on any limit cycle the system is in at least one of the states
$(l_1+d_1,0),$ $(0,l_2+d),$ $(l_1,d),$ $(d,l_2).$
Let us prove that from any of this state the system enters the state of free motion in some time.

Let at the moment of time $ t_0 $ the system is in the state
$$ A(t_0) = (l_1 + d, 0). $$
Then, under the assumptions of the theorem, we have the following sequence of states:
$$A(t_0+d)=(l_1+2d,d),\ A(t_0+n-l_1-d)=(0,n-l_1-d),$$
$$A(t_0+n-l_1)=(d,n-l_1),\ A(t_0+n)=(l_1+d,n).$$
Thus, in $n$ steps, the system returns to the state $(l_1+d, 0),$ i.e. it is in a state of free motion.

Due to symmetry, the state $(0, l_2+d)$ is also a state of free motion.

Suppose that
$$A(t_0)=(l_1,d).$$
So on
$$A(t_0+d-l_1)=(d,2d-l_1),\ A(t_0+n-d)=(n-d+l_1,0),\ A(t_0+n)=(l_1,d).$$

Thus, the state $ (l_1, d_1), $, and by virtue of symmetry, the state $ (d_1, l_2), $ are states of free motion. Theorem~1 is proved.
 \vskip 3pt
 Terms of Theorem 1 are met, for example, for $ n = 10, $ $ l_1 = 1, $ $ l_2 = 2, $ $ d = 3. $
  \ vskip 3pt
In this case, we have a sequence of transitions
 $$(4,0)\to (5,1)\to (6,2)\to (7,3)\to (8,4)\to (9,5)\to $$
 $$\to (0,6)\to (1,7)\to (2,8)\to (3,9)\to (4,0),$$
 $$(0,5)\to (1,6)\to (2,7)\to (3,8)\to (4,9)\to (5,0)\to $$
 $$\to (6,1)\to (7,2)\to (8,3)\to (9,4)\to (0,5),$$
 $$(1,3)\to (2,4)\to (3,5)\to (4,6)\to (5,7)\to ( 6,8)\to $$
 $$\to (7,9)\to (8,0)\to (9,1)\to (0,2)\to (1,3),$$
 $$(3,2)\to (4,3)\to (5,4)\to (6,5)\to (7,6)\to (8,7)\to $$
 $$\to (9,8)\to (0,9)\to (1,0)\to (2,1)\to (3,2),$$

  \subsection*{4.2.System behavior on terms
  $l_2\le d,$ $l_2\le  n-2d,$ $l_1+l_2>n-2d,$}

\hskip 18pt Under terms of the following theorem, the average cluster velocity depends on the initial state.
\ \vskip 3pt
{\bf Theorem 2.} {\it If terms
 $$l_2\le d,\ l_2\le n-2d,\ l_1+l_2>n-2d,$$
 are met then, depending on the initial state of the system, its enters the state of free motion or the average velocity of the clusters is
 $$v=\frac{n}{l_1+l_2+2d}.\eqno(1)$$
 \vskip 3pt
 Proof.} Under the assumptions of Theorem 2, the system, as follows from Lemma~4, cannot be in a state of collapse. If on some limit cycle the system is also not in a state of free motion, then, in accordance with Lemma~3, on any limit cycle the system is in at least one of the states $ (l_1 + d, 0), $ $ (0, l_2 + d ), $ $ (l_1, d), $ $ (d, l_2). $
Let us prove that from any of this state the system enters in some time either the state of free motion, or a limit cycle is realized with the average speed of each cluster, calculated by formula (1).

Suppose that at the moment $ t_0 $ the system is in the state
  $$ A(t_0) = (l_1 + d, 0).$$
  Then we have
 $$A(t_0+n-l_1-d)=(0,n-l_1-d),\ A(t_0+l_2+d)=(0,l_2+d),$$
 $$A(t_0+n)=(n-l_2-d,0),\ A(t_0+l_1+l_2+2d)=(l_1+d,0).$$
   Thus, for a time interval of $l_1+l_2+2d$, all clusters complete a complete revolution, i.e. $n$ movements of each cluster are performed. Thus, on the limit cycle containing the state $l_1+d,$ $l_2+d,$ the average cluster speed is calculated by the formula (1).

 If at the time  $t_0$
 $$A(t_0)=(l_1,d),$$
then we have
 $$A(t_0+d-l_1)=(d,2d-l_1),\ A(t_0+n-l_1)=(0,d-l_1),\ A(t_0+n)=(l_1,d).$$
Therefore, the state $ (l_1, d_2) $ is a state of free movement. Similarly, the $ (d, l_2) $ state is also a free movement state. Thus, each of the four initial states $ (l_1 + d, 0), $ $ (0, l_2 + d), $ $ (l_1, d), $ $ (d, l_2) $ belongs to the limit cycle on which the system either is in a state of free motion, or the average velocity of the clusters is determined by formula (1). Theorem~2 is proved.
\vskip 3pt
  The conditions of Theorem 2 are satistified, for example, for $ n = 12, $ $ l_1 = 2, $ $ l_2 = 4, $ $ d = 4. $

 In this case, we have a sequence of transitions
 $$(6,0)\to (7,1)\to (8,2)\to (9,3)\to (10,4)\to (11,5)\to (0,6)\to (0,7)\to$$
 $$\to (0,8)\to (1,9)\to (2,10)\to (3,11)\to (4,0)\to (5,0)\to (6,0),$$
 $$(2,4)\to (3,5)\to (4,6)\to (5,7)\to (6,8)\to (7,9)\to (8,10)\to $$
 $$\to (9,11) \to (10,0)\to (11,1)\to (0,2)\to (1,3)\to (2,4),$$
 $$(4,4)\to (5,5)\to (6,6) \to (7,7)\to (8,8)\to (9,9)\to (10,10) \to $$
 $$(11,11)\to (0,0)\to (1,1) \to (2,2)\to (3,3)\to (4,4).$$

  \subsection*{4.3. System self-organization on terms
  $l_2\le d,$ $l_1\le n-2d,$ $l_2>  n-2d$}

\hskip 18pt Under terms of the following theorem the system enters the state of self-organization.

\vskip 3pt
 {\bf Theorem 3.} {\it If terms
 $$l_2 \le d,\ l_1\le n-2d,\ l_2>n-2d,$$
 are met then the system  for any initial state is in a state of free motion from some moment in time
 \vskip 3pt
 Proof.} Under the assumptions of Theorem 4, the system, in accordance with Lemma~4, cannot be in a state of collapse. If on some limit cycle the system is also not in a state of free motion, then, in accordance with Lemma~3, on any limit cycle the system is in at least one of the states
  $ (l_1 + d_1,0), $ $ (0, l_2 + d), $ $ (l_1, d), $ $ (d, l_2). $
Let us prove that from any of this state the system enters the state of free motion in some time.

Let at the moment of time $ t_0 $ the system is in the state
$$ A(t_0) = (l_1 + d, 0). $$
Then, under the assumptions of the theorem, we have the following sequence of states:
$$A(t_0+d)=(l_1+2d,d),\ A(t_0+n-l_1-d)=(0,n-l_1-d),$$
$$A(t_0+l_2+d)=(0,l_2+d),\ A(t_0+n)=(n-l_2-d,0),$$
$$A(t_0+l_2+2d)=(d,l_2+2d-n),\ A(t_0+n+l_2)=(d,l_2),$$
$$A(t_0+n+d)=(2d-l_2,d),\ A(t_0+2n+l_2-d)=(0,n-d+l_2),$$
$$A(t_0+2n+l_2)=(d,l_2).$$
So,
$$A(t_0+n+l_2)=A(t_0+2n+l_2),$$
i.e. the system enters a state of free motion.

Thus, the state $ (d, l_2) $ is a state of free motion and from the states $ (l_1 + d, 0), $ $ (0, l_2 + d) $ the system enters the state of free motion.

Let's suppose that
$$A(t_0)=(l_1,d).$$
So on
$$A(t_0+d-l_1)=(d,2d-l_1),\ A(t_0+n-d)=(n+l_1-d,0),$$
$$A(t_0+n-l_1)=(0,d-l_1),\ A(t_0+n)=(l_1,d).$$
Therefore, the state $ (l_1, d) $ is a state of free motion. Theorem~3 is proved.

 \vskip 3pt
 Terms of Theorem 3 is satisfied, for example, for $n=12,$ $l_1=1,$ $l_2=4,$ $d=5.$
 \vskip 3pt
In this case, we have a sequence of transitions
 $$(6,0)\to (7,1)\to (8,2)\to (9,3)\to (10,4)\to (11,5)\to $$
 $$\to (0,6)\to (0,7)\to (0,8)\to (0,9)\to (1,10)\to$$
 $$\to (2,11)\to (3,0)\to (4,1)\to (5,2)\to (5,3)\to (5,4)\to $$
 $$\to (6,5)\to (7,6)\to (8,7)\to (9,8)\to (10,9)\to$$
 $$\to (11,10)\to (0,11)\to (1,0)\to (2,1)\to (3,2)\to$$
 $$\to (4,3)\to (4,5),$$
 $$(1,5)\to (2,6)\to (3,7)\to (4,8)\to (5,9)\to ( 6,10)\to $$
 $$\to (7,11)\to (8,0)\to (9,1)\to (10,2)\to (11,3)\to $$
 $$(0,4)\to (1,5).$$

 \subsection*{4.4. System self-organization on terms
  $l_2\le d,$ $l_1> n-2d$}

\hskip 18pt Under terms of the following theorem the system enters the state of self-organization.

\hskip 18pt
 {\bf Theorem 4.} {\it If terms
 $$l_2 \le d,$$
 $$l_1>n-2d,$$
 are met then the system  for any initial state is in a state of free motion from some moment in time.
 \vskip 3pt
 Proof.} Under the assumptions of Theorem 4, the system, in accordance with Lemma~4, cannot be in a state of collapse. If on some limit cycle the system is also not in a state of free motion, then, in accordance with Lemma~3, on any limit cycle the system is in at least one of the states
  $ (l_1 + d_1,0), $ $ (0, l_2 + d), $ $ (l_1, d), $ $ (d, l_2). $
Let us prove that from any of this state the system enters the state of free motion in some time.

Let at the moment of time $ t_0 $ the system is in the state
$$ A(t_0) = (l_1 + d, 0). $$
Then, under the assumptions of the theorem, we have the following sequence of states:
$$A(t_0+n-l_1-d)=(0,n-l_1-d),\ A(t_0+d)=(l_1+2d-n,d),$$
$$A(t_0+n-d)=(l_1,d),\ A(t_0+n-l_1)=(d,2d-l_1),$$
$$A(t_0+2n-2d)=(n+l_1-d,0),\ A(t_0+2n-l_1-d)=(0,d-l_1),$$
$$A(t_0+2n-d)=(l_1,d).$$
Thus, the state $ (d, l_1) $ is the state of free motion, and from the state $ (l_1 + d, 0) $ the system enters the state of free motion.

Let's suppose that
$$A(t_0)=(0,l_2+d).$$
So on
$$A(t_0+d-l_1)=(d,2d-l_1),\ A(t_0+n-d)=(n-d+l_1,0)$$
$$A(t_0+n-l_1)=(0,d-l_1),\ A(t_0+l_2+d)=(0,l_2+d),$$
$$A(t_0+n)=(l_1,d).$$
Hence, the state $ (l_1, d) $ is the state of free motion
Similarly, the state of free motion is the state $ (d, l_2). $
Theorem~4 is proved.

 \vskip 3pt
 Terms of Theorem 4 are satisfied, for example, for $ n = 12, $ $ l_1 = 3, $ $ l_2 = 4, $ $ d = 5. $
 \vskip 3pt
In this case, we have a sequence of transitions:
 $$(8,0)\to (9,1)\to (10,2)\to (11,3)\to (0,4)\to (1,5)\to $$
 $$\to (2,5)\to (3,5)\to (4,6)\to (5,7)\to (6,8)\to$$
 $$\to (7,9)\to (8,10)\to (9,11)\to (10,0)\to (11,1)\to (0,2)\to $$
 $$\to (1,3)\to (2,4)\to (3,5)\to (4,6)\to (5,7),$$
 $$(0,9)\to (1,10)\to (2,11)\to (3,0)\to (4,1)\to (5,2)\to $$
 $$\to (5,3)\to (5,4)\to (6,5)\to (7,6)\to (8,7)\to$$
 $$\to (9,8)\to (10,9)\to (11,10)\to (0,11)\to (1,0)\to (2,1)\to $$
 $$\to (3,2)\to (4,3)\to (5,4).$$

 \subsection*{4.5. System behavior on terms
  $l_1\le d,$ $l_2>d,$  $l_1+l_2\le  2d$ }

\hskip 18pt Under the conditions of the following theorem, the average cluster velocity depends on the initial state.

\hskip 18pt
 {\bf Theorem 5.} {\it If terms
 $$l_1 \le d,\ l_2>d,\ l_1+l_2\le  2d,$$
  are met then the system  for any initial state is in a state of free motion from some moment in time.
 \vskip 3pt
 Proof.} In accordance with Lemma 4, under the assumptions of Theorem~4, the system cannot be in the state of collapse. If, on a limit cysle, the system is not also in a state of free motion, then, in accordance with Lemma~3, then any limit cycle contains at least one of the states $(l_1+d_1,0),$ $(0,l_2+d),$ $(l_1,d),$ $(d,l_2).$
Let us prove that, from any of these states, the system results in a state of free motion.

Assume that
$$A(t_0)=(l_1+d,0).$$
Then we have
$$A(t_0+d)=(l_1+2d,d),\ A(t_0+n-l_1-d)=(0,n-l_1-d),$$
$$A(t_0+d)=(l_1+2d,d),\ A(t_0+n-l_1-d)=(0,n-l_1-d),$$
$$A(t_0+n-l_1)=(d,n-l_1),\ A(t_0+n-d)=(l_1,n-d),$$
$$A(t_0+n)=(l_1+d,0).$$
Thus the state $(l_1+d,0)$ if a state of free motion.

We can prove analogously that the state $(0,l_2+d)$ is also a state of free motion.

Suppose
$$A(t_0)=(d,l_2).$$
We have
$$A(t_0+n-l_2)=(n-l_2+d,0),\ A(t_0+n-d)=(0,n+l_2-d),$$
$$A(t_0+n)=(d,l_2).$$
Thus the state $(d,l_2)$ is a stare of free movement. Theorem~5 has been proved

 \vskip 3pt
 The condition of Theorem~5 is fulfilled, for example, in the case    $n=12,$ $l_1=1,$ $l_2=3,$ $d=2.$
 \vskip 3pt
We have
 $$(3,0)\to (4,1)\to (5,2)\to (6,3)\to (7,4)\to (8,5)\to $$
 $$\to (9,6)\to (10,7)\to (11,8)\to (0,9)\to $$
 $$\to (1,10)\to (2,11)\to (3,0),$$
 $$(0,5)\to (1,6)\to (2,7)\to (3,8)\to (4,9),$$
 $$(5,10)\to (6,11)\to (7,0)\to (8,1)\to (9,2)\to (10,3)\to $$
 $$\to (11,4)\to (0,5),$$
 $$(1,2)\to (2,3)\to (3,4)\to (4,5)\to (5,6)\to (6,7)\to $$
 $$\to (7,8)\to (8,9)\to (9,10)\to (10,11)\to $$
 $$\to (11,0)\to (0,1)\to (1,2),$$
 $$(2,3)\to (3,4)\to (4,5)\to (5,6)\to (6,7)\to (7,8)\to$$
 $$\to (8,9)\to (9,10)\to (10,11)\to (11,0)\to (0,1)\to$$
 $$\to (1,2)\to (2,3),$$

 \subsection*{4.6. System behavior on terms
  $l_1\le d,$ $d<l_2\le 2d,$  $l_1+l_2\le  2d$ }

\hskip 18pt Under the assumption of the following theorem, the average speed of clusters depends on the initial state.

\hskip 18pt
 {\bf Theorem 6.} {\it Suppose
 $$l_1 \le d,\ d<l_2< 2d,$$
 $$l_1+l_2\le  2d,$$
 
 then $v=1$ or
 $$v_1=v_2=\frac{n}{n-d+l_1+l_2}.\eqno(2)$$
 \vskip 3pt
 Proof.} It follows from Lemma~4 that the system cannot be in the state of collapse. If, on a limit cycle, the system is not in a state of free motion, then, in accordance with Lemma~3, on any limit cycle, the system is in at least of one states $(l_1+d_1,0),$ $(0,l_2+d),$ $(l_1,d),$ $(d,l_2).$ We shall prove that for any of these states the following is true. Either the system results in a state of free motion or the average speed satisfies (2).

Suppose
$$A(t_0)=(l_1+d,0).$$
We have
$$A(t_0+d)=(l_1+2d,d),\ A(t_0+n-l_1-d)=(0,n-l_1-d),$$
$$A(t_0+n-l_1)=(d,n-l_1),\ A(t_0+n)=(l_1+d,0).$$

Thus the state $(l_1+d,0)$ is a state of free motion.

It is proved analogously that the state $(0,l_2+d)$ is a state of free motion.

Assume that
$$A(t_0)=(l_1,d).$$
Then we have
$$A(t_0+d-l_1)=(d,2d-l_1),\ A(t_0+l_2-d)=(d,l_2),$$
$$A(t_0+n-d)=(n+d-l_2,0),\ A(t_0+n+l_2-2d)=(0,l_2-d)$$
$$A(t_0+n)=(2d-l_2,d)\ A(t_0+n-2d+l_1+l_2)=(l_1,d).$$

If
$$A(t_0)=(d,l_2),$$
then
$$A(t_0+n-l_2)=(n-l_2+d,0),\ A(t_0+n-d)=(0,l_2-d),$$
$$A(t_0+n+d-l_2)=(2d-l_2,d)\ A(t_0+n+l_1-d)=(l_1,d).$$
Then, from the state $(d,l_2),$ the system results in the state $(l_1,d),$
which belongs to a limit cycle with the average speed, calculated with  (2). Theorem~6 has been proved.

 \vskip 3pt
 The condition of Theorem 6 is fulfilled if $n=14,$ $l_1=2,$ $l_2=5,$ $d=3.$
 \vskip 3pt
In this case, we have
 $$(5,0)\to (6,1)\to (7,2)\to (8,3)\to (9,4)\to (10,5)\to $$
 $$\to (11,6)\to (12,7)\to (13,8)\to (0,9)\to (1,10)\to $$
 $$\to (2,11)\to (3,12)\to (4,13)\to (5,0),$$
 $$(0,8)\to (1,9)\to (2,10)\to (3,11)\to (4,12)\to (5,13)\to $$
 $$\to (6,0)\to (7,1)\to (8,2)\to (9,3)\to (10,4)\to $$
 $$\to (11,5)\to (12,6)\to (13,7)\to (0,8),$$
 $$(2,3)\to (3,4)\to (3,5)\to (4,6)\to (5,7)\to (6,8)\to $$
 $$\to (7,9)\to (8,10)\to (9,11)\to (10,12)\to (11,13)\to $$
 $$\to (12,0)\to (13,1)\to (0,2)\to (1,3)\to (2,3),$$
 $$(3,5)\to (4,6)\to (5,7)\to (6,8)\to (7,9)\to (8,10)\to $$
 $$\to (9,11)\to (10,12)\to (11,13)\to (12,0)\to (13,1)\to (0,2)$$
 $$\to (1,3)\to (2,3).$$

 \subsection*{4.7. Self-organization in the case $l_1\le d,$ $l_2\ge 2d,$  $_1+l_2\le  2d$}

\hskip 18pt The follows theorem gives a condition of self-organization.

\vskip 3pt
 {\bf Theorem 7.} {\it If
 $$l_1 \le d,\ l_2\ge 2d,\ l_1+l_2\le  2d,$$
 then $v=1.$
\ vskip 3pt
 Proof.} We shall prove that from any of states  $(l_1+d_1,0),$ $(0,l_2+d),$ $(l_1,d),$ $(d,l_2)$ the system results in a state of free motion.

Suppose
$$A(t_0)=(l_1+d,0).$$
Then we have
$$A(t_0+d)=(l_1+2d,d),\ A(t_0+n-l_1-d)=(0,n-l_1-d),$$
$$A(t_0+n-l_1)=(d,n-l_1),\ A(t_0+n)=(l_1+d,0).$$
Therefore the state $(l_1+d,0)$ is a state of free motion.

It is proved analogously that the state $(0,l_2+d)$ is also a state of free motion.

If
$$A(t_0)=(l_1,d),$$
then
$$A(t_0+d-l_1)=(d,2d-l_1),\ A(t_0+l_2-d)=(d,l_2),$$
$$A(t_0+n-d)=(n+d-l_2,0),\ A(t_0+n+l_2-2d)=(0,l_2-d)$$
$$A(t_0+n)=(n+2d-l_2,d),\ A(t_0+n+l_2-2d)=(0,l_2-d),$$
$$A(t_0+n+l_2)=(0,l_2+d).$$
Hence, from the state $(l_1,d),$ the system results in the state    $(0,l_2+d),$ which is a state of free movement.

If
$$A(t_0)=(d,l_2),$$
then
$$A(t_0+n-l_2)=(n-l_2+d,0),\ A(t_0+n-d)=(0,l_2-d),$$
$$A(t_0+n+d)=(0,l_2+d).$$
Therefore, from the state $(d,l_2),$ the system results in the state    $(0,l_2+d),$ which is a state of free movement.

Thus, from the state $(d,l_2),$ the state results in the state in the state $(0,l_2+d),$ which is a state of free motion.
 \vskip 3pt
The conditions of Theorem 7 are satistified, for example, for $n=12,$ $l_1=1,$ $l_2=5,$ $d=2.$
 \vskip 3pt
We have
 $$(3,0)\to (4,1)\to (5,2)\to (6,3)\to (7,4)\to (8,5)\to $$
 $$\to (9,6)\to (10,7)\to (11,8)\to (0,9)\to (1,10)\to $$
 $$\to (2,11)\to (3,0),$$
 $$(0,7)\to (1,8)\to (2,9)\to (3,10)\to (4,11)\to (5,0)\to $$
 $$\to (6,1)\to (7,2)\to (8,3)\to (9,4)\to (10,5)\to $$
 $$\to (11,6)\to (0,7),$$
 $$(1,2)\to (2,3)\to (2,4)\to (2,5)\to (3,6)\to (4,7)\to $$
 $$\to (5,8)\to (6,9)\to (7,10)\to (8,11)\to (9,0)\to $$
 $$\to (10,1)\to (11,2)\to (0,3)\to (1,4)\to (0,5)\to $$
 $$\to (0,6)\to (0,7),$$

 \subsection*{4.8. System behavior on terms
  $l_1\le d,$ $d<l_2< 2d,$  $l_2\le n-2d,$  $l_1+l_2>  n-2d,$  }

\hskip 18pt Under the conditions of the following theorem the system cannot result in a state of free motion or collapse, and there are two possible states.
\vskip 3pt
 {\bf Theorem 8.} {\it If
 $$l_1 \le d,\ d<l_2< 2d,\ l_2\le n-2d,$$
 $$l_1+l_2> n-2d,\eqno(3)$$
 then
 $$v=v_1=v_2=\frac{n}{l_1+l_2+2d}\eqno(4)$$
 or
 $$v=\frac{n}{l_1+l_2+n-2d}.\eqno(5)$$

 \vskip 3pt
 Proof.} It follows from Theorem 8 that the system cannot be in a state of collapse. Then there exist moments such that the leading particle of the cluster~2 is in the cell~0.

 Suppose
 $$A(t_0)=(\alpha,0),$$
and
$$0\le \alpha< d.$$
We have
$$A(d-\alpha)=(d,d-\alpha).$$
Therefore, at time $t=d-\alpha,$ a delay of the cluster~1 begins.

If
$$d\le \alpha< l_1+d,$$
then, at time $t=0,$ the cluster 2 does not move.

If
$$l_1+d<\alpha \le n-d,\eqno(6)$$
then
$$A(n-\alpha)=(0,n-\alpha)$$
and, taking into account  that (3), (6), we get that, at time $t=n-\alpha,$ the cluster~1 does not move.

If
$$\alpha>n-d,$$
then
$$A(d)=(d,d-\alpha)$$
and, therefore, at time $t=d,$ the cluster 1 does not move.

From Lemma 3, it follows that there exists a moment sush that the system is in one of the states $(l_1+d,0),$ $(0,l_2+d),$ $(l_1,d),$ $(d,l_2).$

Suppose
$$A(t_0)=(l_1+d,0).$$
Then we have
$$A(t_0+d)=(l_1+2d,d),\ A(t_0+n-l_1-d)=(0,n-l_1-d),$$
$$A(t_0+l_2+d)=(0,l_2+d),\ A(t_0+l_2+2d)=(d,l_2+2d),$$
$$A(t_0+n)=(n-l_2-d,0),\ A(t_0+l_1+l_2+2d)=(l_1+d,0).$$
Thus the states $(l_1+d,0),$ $(0,l_2+d)$ belong to a limit cycle with the average velocity belongs to a limit cycle with the average velocity satisfied (4).

Assume that
$$A(t_0)=(l_1,d).$$
We have
$$A(t_0+d-l_1)=(d,2d-l_1),\ A(t_0+l_2-d)=(d,l_2),$$
$$A(t_0+n-d)=(n+d-l_2,0),\ A(t_0+n+l_2-2d)=(0,l_2-d),$$
$$A(t_0+n)=(2d-l_2,d),\ A(t_0+n+l_1+l_2-2d)=(l_1,d).$$
Thus the states $(l_1,d),$ $(d,l_2)$ belong to a limit cycle with the average speed satisfied  (5).

 \vskip 3pt
Conditions of Theorem 8 are satisfied if $n=18,$ $l_1=4,$ $l_2=7,$ $d=4.$
 \vskip 3pt
We have
 $$(8,0)\to (9,1)\to (10,2)\to (11,3)\to (12,4)\to (13,5)\to $$
 $$\to (14,6)\to (15,7)\to (16,8)\to (17,9)\to (0,10)\to $$
 $$\to (0,11)\to (1,12)\to (2,13)\to (3,14)\to (4,15)\to$$
 $$\to (5,16)\to (6,17)\to (7,0)\to (8,0),$$
 $$(4,4)\to (4,6)\to (4,7)\to (5,8)\to (6,9)\to (7,10)\to $$
 $$\to (8,11)\to (9,12)\to (10,13)\to (11,14)\to (12,15)\to $$
 $$\to (13,16)\to (14,17)\to (15,0)\to (16,1)\to (17,2)\to $$
 $$\to (0,3)\to (1,4)\to (2,4)\to (3,4)\to (4,4).$$

 \subsection*{4.9.System behavior on terms
  $l_1\le d,$ $2d<l_2\le n-2d,$  $l_1+l_2>  n-2d,$}

\hskip 18pt Under the conditions of Theorem 9, the average velocity does not depend on the initial state.
\vskip 3pt
 {\bf Theorem 9.} {\it If
 $$l_1 \le d,\ 2d\le l_2\le 2n-d,\ l_1+l_2> n-2d,$$
 then there is a single limit cycle with the average velocity that satisfies (4).

 \vskip 3pt
 Proof.} From Lemma 4 it follows that, under the assumptions of Theorem~9, the system does not result in a state of collapse. Therefore there exist moments such that
the leading particle of the cluster~2 is in the cell~0. Assume that
 $$A(t_0)=(\alpha,0).$$

If
$$0\le \alpha< d,$$
then
$$A(d-\alpha)=(d,d-\alpha)$$
 and, therefore, at time $t=d-\alpha,$ a delay of the cluster~1 begins.

If
$$d\le \alpha< l_1+d,$$
then, at time $t=0,$ then the cluster 2 does not move.

If (6) is true, then
$$A(n-\alpha)=(0,n-\alpha).$$
From this and (3), (6), it follows that, at time $t=n-\alpha,$ then the cluster~1 does not move.

If
$$\alpha>n-d,$$
then
$$A(d)=(d,d-\alpha)$$
and, therefore, at time $t=d,$ the cluster~1 does not move.
From this and Lemma~3, it follows that there are exist a moment such that the system is in one of the states  $(l_1+d,0),$ $(0,l_2+d),$ $(l_1,d),$ $(d,l_2).$

If
$$A(t_0)=(l_1+d,0),$$
then
$$A(t_0+d)=(l_1+2d,d),\ A(t_0+n-l_1-d)=(0,n-l_1-d),$$
$$A(t_0+l_2+d)=(0,l_2+d),\ A(t_0+l_2+2d)=(d,l_2+2d),$$
$$A(t_0+n)=(n-l_2-d,0),\ (t_0+l_1+l_2+2d)=(l_1+d,0).$$
Hence the states $(l_1+d,0),$ $(0,l_2+d)$ belong to a limit cycle with the average velocity calculated with the formula (4).

If
$$A(t_0)=(l_1,d),$$
then
$$A(t_0+d-l_1)=(d,2d-l_1),\ A(t_0+l_2-d)=(d,l_2),$$
$$A(t_0+n-d)=(n+d-l_2,0),\ A(t_0+n+l_2-2d)=(0,l_2-d)$$
$$A(t_0+n)=(n+2d-l_2,d),$$
(addition and subtraction modulo $n),$
$$A(t_0+n+l_2-2d)=(0,l_2-d),\ A(t_0+n+l_2)=(0,l_2+d).$$
Thus, from any of states $(l_1,d),$ $(d,l_2),$ the system results in the state $(0,l_2+d),$ which belongs to a limits cycle with the average velocity that satisfies (4).

 \vskip 3pt
 Conditions of Theorem 9 are satisfied if $n=14,$ $l_1=2,$ $l_2=7,$ $d=3.$
 \vskip 3pt
In this case,
 $$(5,0)\to (6,1)\to (7,2)\to (8,3)\to (9,4)\to (10,5)\to $$
 $$\to (11,6)\to (12,7)\to (13,8)\to (0,9)\to (0,10)\to $$
 $$\to (1,11)\to (2,12)\to (3,13)\to (4,0)\to (5,0),$$
 $$(2,3)\to (3,4)\to (3,5)\to (3,6)\to $$
 $$\to (3,7)\to (4,8)\to (5,9)\to (6,10)\to (7,11)\to (8,12)\to $$
 $$\to (9,13)\to (10,0)\to (11,1)\to (12,2)\to (13,3)\to $$
 $$\to (0,4)\to (0,5)\to (0,6)\to (0,7)\to (0,8)\to $$
 $$\to (0,9)\to (0,10).$$

 \subsection*{4.10. System behavior on terms
  $l_1\le d,$ $d<l_2< 2d,$ $l_1\le n-2d,$ $n-2d<l_2\le n-d,$  $l_1+l_2\le  2d$}

 \hskip 18pt  Under the conditions of the following theorem, the system results in a state of free motion from any initial states.

\vskip 3pt
 {\bf Theorem 10.} {\it If
 $$l_1 \le d,\ l_1\le n-2d,\ n-2d<l_2\le n-d,\ l_1+l_2\le  2d,$$
 then the system results in a state of free motion from any initial state.
 \vskip 3pt
 Proof.} From Lemma 4 it follows that, under the assumptions of Theorem~9, the system does not result in a state of collapse. Therefore, in accordance with Lemma~3,  on any limit cycle,  the system is in one of the states  $(l_1+d_1,0),$ $(0,l_2+d),$ $(l_1,d),$ $(d,l_2).$

If
$$A(t_0)=(l_1,d),$$
then
$$A(t_0+d)=(d,2d-l_1),\ A(t_0+n-d)=(n+l_1-d,0),$$
$$A(t_0+n-l_1)=(0,d-l_1),\ A(t_0+n)=(l_1,d).$$
Hence the state $(l_1,d)$ is a state of free motion.

If
$$A(t_0)=(d,l_2),$$
then
$$A(t_0+n-l_2)=(n-l_2+d,0),\ A(t_0+n-d)=(0,l_2-d),$$
$$A(t_0+n+d-l_2)=(2d-l_2,d),\ A(t_0+n)=(d,l_2).$$
Therefore the state $(d,l_2)$ is also a state of free motion.

Suppose
$$A(t_0)=(l_1+d,0).$$
Then we have
$$A(t_0+d)=(l_1+2d,d),\ A(t_0+n-l_1-d)=(0,n+l_1-d),$$
$$A(t_0+l_2+d)=(0,l_2+d),\ A(t_0+l_2+2d)=(d,l_2+2d-n),$$
$$A(n+l_2)=(d,l_2).$$
Thus, from the states $(l_1+d,0),$ $(0,l_2+d),$ the system results in the state $(l_2,d),$ which is a state of free motion.
Theorem~10 has been proved.
 \vskip 3pt
 The condition of Theorem is satisfied if $n=12,$ $l_1=1,$ $l_2=6,$ $d=4.$
 \vskip 3pt
In this, we have the following sequences of states
 $$(1,4)\to (2,5)\to (3,6)\to (4,7)\to (5,8)\to (6,9)\to $$
 $$\to (7,10)\to (8,11)\to (9,0)\to (10,1)\to (11,2)\to $$
 $$\to (0,3)\to (1,4),$$
 $$(4,6)\to (5,7)\to (6,8)\to (7,9)\to (8,10)\to (9,11)\to $$
 $$\to (10,0)\to (11,1)\to (0,2)\to (1,3)\to (2,4)\to $$
 $$\to (3,5)\to (4,6),$$
 $$(5,0)\to (6,1)\to (7,2)\to (8,3)\to (9,4)\to (10,5)\to $$
 $$\to (11,6)\to (0,7)\to (0,8)\to (0,9)\to (0,10)\to $$
 $$\to (1,11)\to (2,0)\to (3,1)\to (4,2)\to (4,3)\to $$
 $$\to (4,4)\to (4,5)\to (4,6).$$

 \subsection*{4.11. System behavior on terms
  $l_1\le d,$ $l_1\le n-2d,$ $d<l_2< 2d,$ $n-2d<l_2\le n-d,$  $l_1+l_2>  2d$}

\hskip 18pt Under the assumptions of Theorem 11, there are two possible values of the average velocity.
\vskip 3pt
 {\bf Theorem 11.} {\it If
 $$l_1 \le d,\ l_1\le n-2d,$$
  $$d<l_2< 2d,\ n-2d<l_2\le n-d,\ l_1+l_2> 2d,$$
 then there exists a single limit cycle, and the average velocity calculatted with the formula (5).

 \vskip 3pt
 Proof.} It follows from Theorem 11 that the system cannot be in a state of collapse. Then there exist moments such that the leading particle of the cluster~2 is in the cell~0.

Assume that at time $t_0$ the system is in the state
 $$(\alpha,0).$$

If
$$0\le \alpha< d,$$
then
$$A(d-\alpha)=(d,d-\alpha).$$
 Therefore, at time $t=d-\alpha,$ a delay of the cluster~1 begins.

If
$$d\le \alpha< l_1+d,$$
then, at time $t=0,$ the cluster 2 does not move.

If
$$l_1+d<\alpha \le n-d,\eqno(7)$$
then
$$A(n-\alpha)=(0,n-\alpha).$$
Ffom this and (3), (7), we get that, at time $t=n-\alpha,$ the cluster~1 does not move.

If
$$\alpha>n-d,$$
then
$$A(d)=(d,d-\alpha)$$
and, therefore, at time $t=d,$ the cluster 1 does not move.

From this and Lemma 3, it follows that there exists a moment such that the system is in one of states $(l_1+d,0),$ $(0,l_2+d),$ $(l_1,d),$ $(d,l_2).$

If
$$A(t_0)=(l_1,d),$$
then
$$A(t_0+d-l_1)=(d,2d-l_1),$$
$$A(t_0+l_2-d)=(d,l_2),$$
$$A(t_0+n-d)=(n+d-l_2,0),$$
$$A(t_0+n+l_2-2d)=(0,l_2-d),$$
$$A(t_0+n)=(2d-l_2,d),$$
$$A(t_0+n+l_1+l_2-2d)=(l_1,d).$$
Therefore the states $(l_1,d),$ $(d,l_2)$ belong to a limit cycle with the average velocity calculated with the formula (5).

If
$$A(t_0)=(l_1+d,0),$$
then
$$A(t_0+d)=(l_1+2d,d),\ A(t_0+n-l_1-d)=(0,n-l_1-d),$$
$$A(t_0+n)=(n-l_2-d,0),\ A(t_0+n+2d-l_2)=(d,2d-l_2).$$
$$A(t_0+l_2)=(d,l_2).$$
Therefore, from the state $(l_1+d,0),$ the system results in a state belonging to a limit cycle with the average velocity calculated with the formula (5).

If
$$A(t_0)=(0,l_2+d),$$
then
$$A(t_0+n-l_2-d)=(n-l_2-d,0),\ A(t_0+d)=(d,l_2+2d-n),$$
$$A(t_0+n-d)=(d,l_2).$$
Thus, from the state $(0,l_2+d),$ the system results in a state belonging to a limit cycle with the average velocity that satisfies  (5).
Theorem~11 has been proved.

 \vskip 3pt
 The condition of Theorem 11 is satisfied if $n=12,$ $l_1=3,$ $l_2=6,$ $d=4.$
 \vskip 3pt
In this case,
 $$(3,4)\to (4,5)\to (4,6)\to (5,7)\to (6,8)\to (7,9)\to $$
 $$\to (8,10)\to (9,11)\to (10,0\to (11,1)\to (0,2)\to $$
 $$\to (1,3)\to (2,4)\to (3,4),$$
 $$\to (5,16)\to (6,17)\to (7,0)\to (8,0),$$
 $$(7,0)\to (8,1)\to (9,2)\to (10,3)\to (11,4)\to (0,5)\to $$
 $$\to (0,6)\to (0,7)\to (0,8)\to (0,9)\to (0,10)\to $$
 $$\to (1,11)\to (2,0)\to (4,2)\to (4,3)\to (4,4)\to $$
 $$\to (4,5)\to (4,6),$$
 $$(0,10)\to (1,11)\to (2,0)\to (3,1)\to (4,2)\to (4,3)\to $$
 $$\to (4,4)\to (4,5)\to (4,6).$$

\subsection*{4.12. System behavior on terms
  $l_1\le d_2,$ $l\leq n-2d$, $l_2> n-d,$ $d<l_2< 2d$ }

\hskip 18pt Under the conditions of the following theorem, the system results in a state of free motion from any initial state.
\vskip 3pt
 {\bf Theorem 12.} {\it Suppose
 $$l_1 \le d,\ l_1\le n-2d,\ l_2> n-d,$$
  $$d<l_2< 2d,\ l_1+l_2\le 2d.$$
 Then the system results in the state oof free motion from any initial state.

\vskip 3pt
 Proof.} From Lemma 4 it follows that, under the assumptions of Theorem~9, the system does not result in a state of collapse. Therefore, in accordance with Lemma~3,  on any limit cycle,  the system is in one of the states  $(l_1+d_1,0),$ $(0,l_2+d),$ $(l_1,d),$ $(d,l_2).$ We shall prove that the system results in a state of free motion from any initial state.

If
$$A(t_0)=(l_1,d),$$
then
$$A(t_0+d-l_1)=(d,2d-l_1),\ A(t_0+n-d)=(n-d+l_1,0),$$
$$A(t_0+n-l_1)=(0,d-l_1),\ A(t_0+n-l_1)=(l_1,d).$$
Whence the state $(l_1,d)$ is a state of free motion.

Suppose
$$A(t_0)=(d,l_2).$$
Then
$$A(t_0+n-l_2)=(n-l_2+d,0),\ A(t_0+n-d)=(0,l_2-d),$$
$$A(t_0+n+d-l_2)=(2d-l_2,d),\ A(t_0+n)=(d,l_2).$$
Therefore the state $(l,d_2)$ is also the state of free motion.

Assume that
$$A(t_0)=(0,l_2+d-n).$$
Then
$$A(t_0+n-l_2-d)=(n-l_2-d,0),\ A(t_0+d)=(d,l_2+2d-n),$$
$$A(t_0+n-d)=(d,l_2).$$
Therefore, from the state $(0,l_2+d),$ the system results in the state $(d,l_2),$ which is a state of free motion.

If
$$A(t_0)=(l_1+d,0),$$
then
$$A(t_0+n-l_2)=(n-l_2,d),\ A(t_0+d)=(d,l_2+2d-n),$$
$$A(t_0+n-d)=(d,l_2).$$
Thus, from the state $(l_1+d,0),$ the system results in a state $(0,l_2+d),$ from which the system results in the state $(d,l_2),$ and the state $(d,l_2)$ is a state of free movement. Theorem~12 has been proved.

 \vskip 3pt
 The condition of Theorem 12 is satisfied for $n=12,$ $l_1=1,$ $l_2=9,$ $d=5.$
 \vskip 3pt
In this case,
$$(1,5)\to (2,6)\to (3,7)\to (4,8)\to (5,9)\to (6,10)\to $$
 $$\to (7,11)\to (8,0)\to (9,1)\to (10,2)\to (11,3)\to $$
 $$\to (0,4)\to (1,5),$$
 $$(6,0)\to (7,1)\to (8,2)\to (9,3)\to (10,4)\to (11,5)\to $$
 $$\to (0,6)\to (0,7)\to (0,8)\to (0,9)\to (0,10)\to $$
 $$(0,11)\to (0,0)\to (0,1)\to (0,2)\to (1,3)\to (2,4)\to $$
 $$\to (3,5)\to (4,6)\to (5,7)\to (5,8)\to (5,9)\to $$
 $$\to (6,10)\to (7,11)\to (8,0)\to (9,1)\to (10,2)\to $$
 $$\to (11,3)\to (0,4)\to (1,5)\to (2,6)\to (3,7)\to $$
  $$\to (4,8)\to (5,9).$$

 \subsection*{4.13. System behavior on terms
  $l_1\le d,$ $l_1\le n-2d,$ $d<l_2< 2d,$ $l_2>n-d,$  $2d <l_1+l_2\le n$}

\hskip 18pt Under the conditions of the following theorem, the system results in a state of free motion from any initial state.
\vskip 3pt
 {\bf Theorem 13.} {\it Assume that
 $$l_1 \le d,\ l_1\le n-2d,\ l_2> n-d,$$
  $$d<l_2< 2d,\  l_2> n-d,\ l_1+l_2> 2d.$$
 Then the system results in the state of free motion from any iinitial state.
\vskip 3pt
 Proof.} From Lemma 4 it follows that, under the assumptions of Theorem~9, the system does not result in a state of collapse. Therefore, in accordance with Lemma~3,  on any limit cycle,  the system is in one of the states  $(l_1+d_1,0),$ $(0,l_2+d),$ $(l_1,d),$ $(d,l_2).$ We shall prove that the system results in a state of free motion from any initial state.

Suppose
$$A(t_0)=(l_1,d).$$
We have
$$A(t_0+d-l_1)=(d,2d-l_1),\ A(t_0+n+l_2-d)=(d,l_2),\ A(t_0+n-d)=(n-l_2+d,0),$$
$$A(t_0+n+l_2-2d)=(0,l_2-d),\ A(t_0+n+l_2-d)=(d,l_2).$$
Hence the state $(d,l_2)$ is a state of free motion, and, from the state  $(l_1,d),$ the system results in the state of free motion.

$$A(t_0)=(l_1+d,0).$$
Then we have the following sequence of states
$$A(t_0+n+d)=(l_1+2d,d),\ A(t_0+n-l_1-d)=(0,n-l_1-d),\ A(t_0+n+l_2+d)=(0,l_2+d-n),$$
$$A(t_0+n+d)=(n-l_2,d),\ A(t_0+l_2+2d)=(d,l_2+2d-n),\ A(l_2+d+n)=(0,l_2+d-n).$$
Thus the state $(0,l_2+d-n)$ if a state of free motion, and, from the state $(l_1+d,d),$ the system results in a state of free motion.

If
$$A(t_0)=(0,l_2+d-n),$$
then
$$A(t_0+n-l_2-d)=(n-l_2-d,0),\ A(t_0+d)=(d,l_2+2d-n),\ A(t_0+n-d)=(d,l_2).$$
Thus, from the state $(0,l_2+d),$ the system results in the state $(d,l_2),$ which is a state of free motion.
\vskip 3pt
The condition of Theorem 13 is satisfied if $n=20,$ $l_1=3,$ $l_2=14,$ $d=8.$
\vskip 3pt
We have
 $$(11,0)\to (12,1)\to (13,2)\to (14,3)\to (15,4)\to (16,5)\to $$
 $$\to (17,6)\to (18,7)\to (19,8)\to (0,9)\to (0,10)\to $$
 $$\to (0,11)\to (0,12)\to (0,13)\to (0,14)\to (0,15)\to (0,16)\to $$
 $$\to (0,17)\to (0,18)\to (0,19)\to (0,0)\to (0,1)\to $$
 $$(0,2)\to (1,3)\to (2,4)\to (3,5)\to (4,6)\to (5,7)\to $$
 $$\to (6,8)\to (7,9)\to (8,10)\to (9,11)\to (10,12)\to $$
 $$\to (11,13)\to (12,14)\to (13,15)\to (14,16)\to (15,17)\to $$
 $$\to (16,18)\to (17,19)\to (18,0)\to (19,1)\to (0,2),$$
  $$(3,8)\to (4,9)\to (5,10)\to (6,11)\to (7,12)\to (8,13)\to $$
 $$\to (8,14)\to (9,15)\to (10,16)\to (11,17)\to (12,18)\to $$
 $$\to (13,19)\to (14,0)\to (15,1)\to (16,2)\to (17,3)\to (18,4)\to $$
 $$\to (19,5)\to (0,6)\to (1,7)\to (2,8)\to (3,9)\to $$
 $$(4,10)\to (5,11)\to (6,12)\to (7,13)\to (8,14). $$

  \subsection*{4.14. System behavior on terms
  $l_1\le d,$ $l_1\le n-2d,$ $l_2>n-d,$  $l_1+l_2> n,$}

 \hskip 18pt Under the condition of the following theorem, there is a single limit cycle, and the value of the cluster~1 velocity is equal to the value of the cluster~2 multiplied by~2.
\vskip 3pt
 {\bf Theorem 14.} {\it If
 $$l_1 \le d,\ l_1\le n-2d,\ l_2> n-d,\ l_1+l_2> n,$$
 then the velocities of the clusters are equal to
 $$v_1=\frac{n}{2(l_1+l_2)},\eqno(8)$$
 $$v_2=\frac{n}{l_1+l_2}.\eqno(9)$$

\vskip 3pt
Proof.} In accordance with Lemmas~1 and 4, the system cannot be in a state of collapse or free motion. In accordance with Lemma~3, on a limit cycle results in one of the states  $(l_1+d_1,0),$ $(0,l_2+d-n),$ $(l_1,d),$ $(d,l_2).$ We shall prove that average velocities of clusters are calculated with (8), (9).

If
$$A(t_0)=(l_1+d,0),$$
then we have
$$A(t_0+d)=(l_1+2d,d),\ A(t_0+n-l_1-d)=(0,n-l_1-d),$$
$$A(t_0+n)=(0,0),\ A(t_0+l_2+d)=(0,l_2+d-n),$$
$$A(t_0+n+d)=(n-l_2,d),\ A(t_0+l_1+l_2+d)=(l_1,d),$$
$$A(t_0+l_2+2d)=(d,2d-l_1),\ A(t_0+l_1+2l_2)=(d,l_2),$$
$$A(t_0+n+l_1+l_2)=(n+d-l_2,0),\ A(t_0+2(l_1+l_2))=(l_1+d,0).$$
Thus there exists a single limit cycle. The period of the cycle equals $2(l_1+l_2).$ On the cycle, there are two turns of the cluster~1 and a turn of the cluater~2. From this, Theorem~14 follows.
\vskip 3pt
For example, the condition of Theorem 14 is satisfied if $n=12,$ $l_1=2,$ $l_2=11,$ $d=3.$
  \vskip 3pt
In this case, we have the following sequence of transitions
 $$(5,0)\to (6,1)\to (7,2)\to (8,3)\to (9,4)\to (10,5)\to $$
 $$\to (11,6)\to (0,7)\to (0,8)\to (0,9)\to (0,10)\to $$
 $$\to (0,11)\to (0,0)\to (0,1)\to (0,2)\to (1,3)\to $$
 $$\to (2,3)\to (3,4)\to (3,5)\to (3,6)\to (3,7)\to (3,8)\to  $$
 $$\to (3,9)\to (3,10)\to (3,11)\to (4,0)\to (5,0).$$

 \subsection*{4.15. System behavior on terms
$l_1 \le d,$  $n-2d<l_1\le n-d,$ $l_2>n-d,$ $l_1+l_2>n$}

 \hskip 18pt If the condition of the following theorem there exists a single limit cycle with the average velocity less than~1.
 \vskip 3pt
{\bf Theorem 15.} {\it If
  $$l_1\le d,\  n-2d<l_1<n-d,\    l_2>n-d,$$
 $$l_1+l_2>n,$$
 then, from any initial state, the system results in a state of a single limit cycle with the average velocities of clusters equal to
 $$v_1=v_2=\frac{n}{n-2d+l_1+l_2}.\eqno(9)$$
 \vskip 3pt
  Proof.} Under the assumptions of Theorem 15, in accordance with Lemma~1, the system cannot be in a state of free motion and,  in accordance with Lemma~1, the system cannot be in a state of

  If
  $$A(t_0)=(l_1+d,0),$$
  then
  $$A(t_0+n-l_1-d)=(0,n-l_1-d),\ A(t_0+d)=(l_1+2d-n,d),$$
 $$A(t_0+n-l_1)=(d,n-l_1),\ A(t_0+l_2)=(d,l_2),$$
 $$A(t_0+n)=(n+d-l_2,0),\ (t_0+n+l_2-d)=(0,l_2-d),$$
 $$A(n+d)=(2d-l_2,d),\ A(n+l_1+l_2-d)=(l_1,d),$$
$$A(n+l_2)=(d,2d-l_1),\ A(n+l_1+2l_2-2d)=(d,l_2).$$
Therefore the states $(l_1,d),$ $(d,l_2)$ belong to a limit cycle such that the average velocity is calculated with the formula (9).

 If
  $$A(t_0)=(0,l_2+d-n),$$
  then
  $$A(t_0+d)=(d,l_2+2d-n),\ A(t_0+n-d)=(d,l_2).$$
 Thus, from the state $(0,l_2+d-n),$ the system results in the state, which belongs to a limit cycle with the average velocity satisfied~(9). Theorem~15 has been proved.

 \subsection*{4.16. System behavior on terms
$l_1 \le d,$  $l_1\le n-2d,$ $n-2d <l_2\le n-d,$ $l_2\ge 2d,$ $l_1+l_2\le n$}

  \hskip 18pt If the conditions of the following theorem holds, then there exists a single limit cycle such that the average velocity of the cluster~1 is equal to~1/2, and the velocityaverage velocity of the cluster~2  with speed 1.
 \vskip 3pt
{\bf Theorem 16.} {\it If
  $$l_1\le d,\  l_1\le n-2d,\ l_2> n-d,\ l_2\ge 2d,\  l_1+l_2\le n,$$
 then, from any initial state, the system results in a state such that this state belongs to a limit cycle with the average velicities of clusters
 $$v_1=\frac{1}{2},\ v_2=1.\eqno(10)$$
 \vskip 3pt
 Proof.} Under the assumptions of Lemma~1, the system cannot be in a state of free movement, and, in accordance with Lemma~4, the system cannot be in a stete of collapse.

 If
  $$A(t_0)=(l_1+d,0),$$
 In this case,
  $$A(t_0+n-l_1-d)=(0,n-l_1-d),\  A(t_0+n)=(0,0),$$
 $$A(t_0+l_2+d)=(0,l_2+d-n),\ A(t_0+n+d)=(n-l_2,d),$$
 $$A(t_0+l_2+2d)=(d,l_2+2d-n),\ A(t_0+n+l_2)=(d,l_2),$$
 $$A(t_0+2n)=(n+d-l_2,0),\ A(2n+d)=(n+2d-l_2,d),$$
 $$A(2n+l_2-d)=(0,l_2-d),\ A(2n+l_2+d)=(0,l_2+d-n).$$
Whence the states $(l_1+d,0),$ $(0,l_2+2d-n),$ $(d,l_2)$ belong to a limit cycle with the average velocity satisfied (10).

  Soppuse
  $$A(t_0)=(l_1,d)$$
  Then we have
  $$A(t_0+d-l_1)=(d,2d-l_1),\ A(t_0+l_2)=(d,l_2).$$
 Thus, from the state $(l_1,d),$ the system results in the state $(d,l_2),$ which belongs to a limit cycle with the average veelocities of clusters, claculated with the formula (10). Theorem~16 has been proved.

\vskip 3pt

For example, the condition of Theorem~16 is satisfied for $n=12,$ $l_1=1,$ $l_2=10,$ $d=3.$

We have
$$(4,0)\to (5,1)\to (6,2)\to (7,3)\to (8,4)\to (9,5)\to $$
$$\to (10,6)\to (11,7)\to (0,8)\to (0,9)\to (0,10)\to (0,11)\to $$
$$\to (0,0)\to (0,1)\to (1,2)\to (2,3)\to (3,4)\to (3,5)\to $$
$$\to (3,6)\to (3,7)\to (3,8)\to (3,9)\to (3,10)\to (4,11)\to $$
$$\to (5,0)\to (6,1)\to (7,2)\to (8,3)\to (9,4)\to (10,5)\to$$
$$\to (11,6)\to (0,7)\to (0,8\to (0,9)\to (0,10)\to (0,11)\to $$
$$\to (0,0)\to (0,1),$$
$$(1,3)\to (2,4)\to (3,5)\to (4,6)\to (5,7)\to (5,8)\to $$
$$\to (5,9)\to (5.10.$$

  \subsection*{4.17. System behavior on terms
  $d<l_1< 2d,$ $d<l_2< 2d,$ $l_1+l_2\le n-2d$}

\hskip 18pt If th condition of the following theorem holds, the average velocities of both clusters are equal to  0 or 1.
\vskip 3pt
{\bf Theorem 17.} {\it Suppose
 $$d<l_1<2d,\ d<l_2<2d,$$
 $$l_1+l_2\le n-d_1-d_2,$$
then, depending on the initial state, the system results in the state of free motiion or collapse.
 \vskip 3pt
 Proof.} If, on a limit cycle, the system is not in the state neither free movement or collapse, then, in accordance with Lemma~3, on any limit cycle, the system results in at least on of the states
 $(l_1+d,0),$ $(0,l_2+d),$ $(l_1,d),$ $(d,l_2).$
We shall prove that, from this state, the system results in the state of free motion or collapse.

 If
  $$A(t_0)=(l_1+d,0),$$
 then we have
 $$A(t_0+d)=(l_1+2d,d),\ A(t_0+n-l_1-d)=(0,n-l_1-d),$$
 $$A(t_0+n-l_1)=(d,n-l_1),\ A(t_0+n)=(l_1+d,0).$$
 Thus the state $(l_1+d, 0)$ is a state of free motion.

 Assume that
 $$A(t_0)=(l_1,d).$$
Then we have the following sequences of states
 $$A(t_0+n-l_1)=(0,n-l_1+d),\ A(t_0+n-l_1+d)=(d,2d-l_1),$$
 $$A(t_0+n)=(d,d).$$
Therefore, from the state $(l_1,d),$ the system results in the state $(d,d),$
which is a state of collapse.

 Assume that
 $$A(t_0)=(d,l_2).$$
 Then the following sequence of states is realized
 $$A(t_0+n-l_2)=(n-l_2+d,0),\ A(t_0+n-l_2+d)=(2d-l_2,d),$$
 $$A(t_0+n)=(d,d).$$
 Thus, from the state $(l_1,d),$ the state results in the state $(d,d),$ i.e., the system results in the state of collapse. Theorem~17 has been proved. Theorem~17 has been proved.
 \vskip 3pt
 The condition of Theorem 17 is satisfied, for example, in the case $n=12,$ $l_1=3,$ $l_2=3,$ $d=2.$

 We have the following sequence of transitions
 $$(5,0)\to (6,1)\to (7,2)\to (8,3)\to (9,4)\to (10,5)\to $$
 $$\to (11,6)\to (0,7)\to (1,8)\to (2,9)\to (3,10)\to $$
 $$\to (4,11)\to (5,0),$$
 $$(3,2)\to (4,3)\to (5,4)\to (6,5)\to (7,6)\to (8,7)\to $$
 $$\to (9,8)\to (10,9)\to (11,10)\to (0,11)\to (1,0)\to (2,1)\to $$
 $$\to (3,2).$$

  \subsection*{4.18.System behavior on terms
  $d<l_1< 2d,$ $l_2\ge 2d,$ $l_1+l_2\le n-2d$}

\hskip 18pt Under the condition of the following theorem, depending on the
initial state, the average velocity of clusters is equal to 0 or 1.
\vskip 3pt
{\bf Theorem 18.} {\it If the conditions
 $$d<l_1\le d,\ l_2\ge 2d,\ l_1+l_2\le n-d_1-d_2$$
 holds, then, depending on the initial state, the system results in a state
 of free movement or collapse.
 \vskip 3pt
 Proof.} If, on a limit cycle, the system is not in a state of free motion or collapse, then, in accordance with Lemma~3, on this limit cycle, the system results in at least one of states $(l_1+d,0),$ $(0,l_2+d),$ $(l_1,d),$ $(d,l_2).$
We shall prove that, from any of these state, the system results in a state of free motion or collapse.

 Suppose
  $$A(t_0)=(l_1+d,0).$$
Then we have
 $$A(t_0+d)=(l_1+2d,d),\ A(t_0+n-l_1-d)=(0,n-l_1-d),$$
 $$A(t_0+n-l_1)=(d,n-l_1),\ A(t_0+n)=(l_1+d,0).$$
 Whence the state $(l_1+d, 0)$ is a state of free motion.

 Assume that
 $$A(t_0)=(l_1,d).$$
We have
 $$A(t_0+n-l_1)=(0,n-l_1+d),\ A(t_0+n-l_1+d)=(d,2d-l_1),$$
 $$A(t_0+n)=(d,d).$$
Whence, from the state $(l_1,d),$ the systems results in a state of free $(d,d),$ which is a state of free motion.

If
 $$A(t_0)=(d,l_2),$$
then we have
 $$A(t_0+n-l_2)=(n-l_2+d,0),\ A(t_0+n)=(0,l_2-d),$$
 $$A(t_0+n+2d)=(0,l_2+d).$$
Hence, from the state $(d,l_2),$ the system results in the state $(0,l_2+d),$
which is the state of free motion.

The state $(d,d)$ is a state of collapse. Theorem~18 has been proved.
 \vskip 3pt
The condition of Theorem 18 holds if, for example, $n=12,$ $l_1=3,$ $l_2=4,$ $d=2.$

 In this case, we have
 $$(5,0)\to (6,1)\to (7,2)\to (8,3)\to (9,4)\to (10,5)\to $$
 $$\to (11,6)\to (0,7)\to (1,8)\to (2,9)\to (3,10)\to $$
 $$\to (4,11)\to (5,0),$$
 $$(0,6)\to (1,7)\to (2,8)\to (3,9)\to (4,10)\to (5,11)\to $$
 $$\to (6,0)\to (7,1)\to (8,2)\to (9,3)\to (10,4)\to (11,5)\to $$
 $$\to (0,6),$$
 $$(3,2)\to (4,3)\to (5,4)\to (6,5)\to (7,6)\to (8,7)\to $$
 $$\to (9,8)\to (10,9)\to (11,.10)\to (0,11)\to (1,0)\to $$
 $$\to (2,1)\to (2,2).$$
 $$(2,4)\to (3,5)\to (4,6)\to (5,7)\to (6,8)\to (7,9)\to $$
 $$\to (8,10)\to (9,11)\to (10,0)\to (11,1)\to (0,2)\to $$
 $$\to (0,3)\to (0,4)\to (0,5)\to (0,6).$$

  \subsection*{4.19. System behavior on terms
  $l_1\ge 2d,$ $l_2\ge 2d,$ $l_1+l_2\le n-2d$}

\hskip 18pt Under the assumptions of the following theorem, depending on the initial state, the average of clusters is equal to 0 or 1.
\vskip 3pt
{\bf Theorem 19.} {\it If
 $$l_2\ge 2d,$$
 $$l_1+l_2\le n-2d$$
 then, depending on the initial state, the system results in the state of free motion or collapse.
 \vskip 3pt
 Proof.} If, on a limit cycle, the system is not in a state of free motion or collapse, then, in accordance with Lemma~3, this system results in one of states $(l_1+d,0),$ $(0,l_2+d),$ $(l_1,d),$ $(d,l_2).$ We shall prove that, from any of these states, the system results in a state of free motion or collapse.

 Assume that
  $$A(t_0)=(l_1+d,0),$$
 Then we have the following sequence of transitions
 $$A(t_0+d)=(l_1+2d,d),$$
 $$A(t_0+n-l_1-d)=(0,n-l_1-d),$$
 $$A(t_0+n-l_1)=(d,n-l_1),$$
 $$A(t_0+n)=(l_1+d,0).$$
 Therefore the state $(l_1+d, 0)$ is a state of free motion.

 Suppose
  $$A(t_0)=(0,l_2+d).$$
 Then we have
 $$A(t_0+d)=(0,l_2+2d),$$
 $$A(t_0+n-l_2-d)=(n-l_2-d,0),$$
 $$A(t_0+n-l_2)=(n-l_2,d),$$
 $$A(t_0+n)=(l_2+d,0).$$
 Hence the state $(0,l_2+d)$ is a state of free motion.

Suppose
 $$A(t_0)=(l_1,d).$$
 Then we have the following sequence of states
 $$A(t_0+n-l_1)=(0,n-l_1+d),$$
 $$A(t_0+n)=(l_1-d,0),$$
 $$A(t_0+n+2d)=(l_1+d,0).$$
 Therefore, from the state $(l_1,d),$ the system results in the state $(l_1+d,d),$ which is a state of free movement.

 Suppose
 $$A(t_0)=(d,l_2).$$
 Then we have
 $$A(t_0+n-l_2)=(n-l_2+d,0),$$
 $$A(t_0+n)=(0,l_2-d),$$
 $$A(t_0+n+2d)=(0,l_2+d).$$
Thus, from the state $(d,l_2),$ the system results in the state $(0,l_2+d),$ which is a state of free movement. Theorem~19 has been proved.

 \vskip 3pt
 The condition of Theorem 19 holds on, for example, when $n=12,$ $l_1=6,$ $l_2=7,$ $d=2.$

 We have the following sequence of states
 $$(8,0)\to (9,1)\to (10,2)\to (11,3)\to (0,4)\to (1,5)\to $$
 $$\to (2,6)\to (3,7)\to (4,8)\to (5,9)\to (6,10)\to $$
 $$\to (7,11)\to (5,0),$$
 $$(0,6)\to (1,7)\to (2,8)\to (3,9)\to (4,10)\to (5,11)\to $$
 $$\to (6,0)\to (7,1)\to (8,2)\to (9,3)\to (10,4)\to (11,5)\to $$
 $$\to (0,6),$$
 $$(3,2)\to (4,3)\to (5,4)\to (6,5)\to (7,6)\to (8,7)\to $$
 $$\to (9,8)\to (10,9)\to (11,.10)\to (0,11)\to (1,0)\to $$
 $$\to (2,1)\to (2,2).$$
 $$(2,4)\to (3,5)\to (4,6)\to (5,7)\to (6,8)\to (7,9)\to $$
 $$\to (8,10)\to (9,11)\to (10,0)\to (11,1)\to (0,2)\to $$
 $$\to (0,3)\to (0,4)\to (0,5)\to (0,6).$$

   \subsection*{4.20. System behavior on terms
  $l_1\ge 2d,$  $l_2\le n-2d,$
  $n-2d<l_1+l_2\le n,$
  }
\hskip 18pt Under the assuptions of the following theorem, depending on the initial state, the velocity of both clusters is equal to 0 or $v=\frac{n}{l_1+l+2+2d}.$
\vskip 3pt
{\bf Theorem 20.} {\it If the conditions
 $$l_1\ge 2d,\  l_2\le n-2d,$$
 $$n-2d<l_1+l_2\le n,$$
 holds, then, depending on the initial state, either the system results in the state of collapse or a limit cycle with the  average velocity calculted with (1) is realized.
 \vskip 3pt
Proof.} In accordance with Lemma~1, if the condition of Theorem  20, then, from any initial state, the system does not results in a sraye of free movement. If, on a limit cycle, the system is not in a state of free motion or collapse, then, in accordance with Lemma~3, this system results in one of states $(l_1+d,0),$ $(0,l_2+d),$ $(l_1,d),$ $(d,l_2).$ We shall prove that, from any of these states, the either the system results in a state of free motion or the average velocity of clusters satisfies (1).

 If
  $$A(t_0)=(l_1+d,0),$$
 then we have
 $$A(t_0+d)=(l_1+2d,d),\ A(t_0+n-l_1-d)=(0,n-l_1-d),$$
 $$A(t_0+l_2+d)=(0,l_2+d),\  A(t_0+l_1+l_2+2d)=(d_1+l,0).$$
 Therefore, the states $(l_1+d_1, 0)$ and $(0,l_2+d_2)$ belong to a limit cycle with the average velocity calculated with the formula (1).

 If
 $$A(t_0)=(l_1,d),$$
 then we have
 $$A(t_0+n-l_1)=(0,n-l_1+d),\ A(t_0+l_2+d)=(0,l_2+d).$$
 Therefore, from the state $(l_1,d_2),$ the system results in a state $(0,l_2+d),$ and, from this state, the system results in a state pf a limit cycle with the average velocity calculated with the formula (1). From the state $(d,l_2),$ the system results in a state $(l_1+d,0),$ from which the system also results in a state of a limit cycle calculated with the  with the average velocity  formula (1).

 The state  $(d,d)$ is a state of collapse. Theorem~ 20 has been proved.
\vskip 3pt
 The condition of Theorem 20 holds, for example, in the case $n=12,$ $l_1=4,$ $l_2=5,$ $d=2.$

 We have the following sequences of transitions
  $$(6,0)\to (7,1)\to (8,2)\to (9,3)\to (10,4)\to (11,5)\to $$
 $$\to (0,6)\to (0,7)\to (1,8)\to (2,9)\to (3,10)\to $$
 $$\to (4,11)\to (5,0)\to (6,0),$$
 $$(4,2)\to (5,3)\to (6,4)\to (7,5)\to (8,6)\to (9,7)\to $$
 $$\to (10,8)\to (11,9)\to (0,10)\to (1,11)\to (2,0)\to (11,5)\to $$
 $$\to (0,6)$$
 $$(3,2)\to (4,3)\to (5,4)\to (6,5)\to (7,6)\to (8,7)\to $$
 $$\to (9,8)\to (10,9)\to (11.10)\to (0,11)\to (1,0)\to $$
 $$\to (2,1)\to (3,2).$$
 $$(2,4)\to (3,5)\to (4,6)\to (5,7)\to (6,8)\to (7,9)\to $$
 $$\to (8,10)\to (9,11)\to (10.0)\to (11,1)\to (0,2)\to $$
 $$\to (1,3)\to (2,4).$$

 \subsection*{4.21. Behavior of the systems in the case
  $l_1\ge 2d,$ $l_2\le n-2d,$
  $l_1+l_2> n,$
  }
\hskip 18pt  If the condition of the following theorem holds, then, depending on the initial state, the average velocity equals 0 or
$v=\frac{n}{l_1+l+2+2d}.$
\vskip 3pt
{\bf Theorem 21.} {\it If the conditions
 $$l_1\ge 2d,\  l_2\le n-2d,\ l_1+l_2> n,$$
 holds, then, depending on the initial state, either the system results in the state of collapse or a limit cycle with the  average velocity calculted with (1) is realized.
  \vskip 3pt
 Proof.} In accordance with Lemma~1, if the condition of Theorem  20, then, from any initial state, the system does not results in a sraye of free movement. If, on a limit cycle, the system is not in a state of free motion or collapse, then, in accordance with Lemma~3, this system results in one of states $(l_1+d,0),$ $(0,l_2+d),$ $(l_1,d),$ $(d,l_2).$ We shall prove that, from any of these states, the either the system results in a state of free motion or the average velocity of clusters satisfies (1).

 Suppose
  $$A(t_0)=(l_1+d,0).$$
  Then we have
  $$A(t_0+d)=(l_1+2d,0),$$
  $$A(t_0+n-l_1-d)=(0,n-l_1-d),$$
 $$A(t_0+l_2+d)=(0,l_2+d),$$
 $$A(t_0+l_2+2d)=(d,l_2+2d),$$
  $$A(t_0+n)=(n-l_2-2d,0),$$
  $$A(t_0+l_1+l_2+2d)=(l_1+d,0).$$
 Therefore the states $(l_1+d, 0)$ and $(0,l_2+d)$ belong to a limit cycle with the average velocity that is calculated, which is calculated with (1).

 Suppose
 $$A(t_0)=(l_1,d).$$
 Then we have
 $$A(t_0+n-l_1)=(0,n-l_1+d),$$
 $$A(t_0+l_2+d)=(0,l_2+d).$$

Thus, from the state $(l_1,d),$ the system results in the state $(0,l_2+d),$ from which the system results in a state of a limit cycle with the average velocity calculted with (1). From the state $(d,l_2)$ the system results in the state $(l_1+d,0),$ from which the system results in a state of a limit cycle with the average velocity calculted with (1).

 The state $(d,d)$ is a state of collapse. Theorem 21 has been proved.

\vskip 3pt
 The condition of Theorem 21 holds for $n=12,$ $l_1=6,$ $l_2=7,$ $d=2.$

 In this case, the following sequences of the transitions are realized
 $$(8,0)\to (9,1)\to (10,2)\to (11,3)\to (0,4)\to (0,5)\to $$
 $$\to (0,6)\to (0,7)\to (0,8)\to (0,9)\to (1,10)\to $$
 $$\to (2,11)\to (3,0)\to (4,0)\to (5,0)\to (6,0)\to (7,0\to (8,0),$$
 $$(6,2)\to (7,3)\to (8,4)\to (9,5)\to (10,6)\to (11,7)\to $$
 $$\to (0,8)\to (0,9),$$
 $$(2,7)\to (3,8)\to (4,9)\to (5,10)\to (6,11)\to (7,0)\to $$
 $$\to (8,0).$$

 \subsection*{4.22. Resulting in a state of collapse in the case
  $l_2\ge 2d,$ $l_1\le n-2d,$ $n-2d<  l_2\le n-d$
  }
\hskip 18pt Under the assumptions of the following theorem, the system results in a state of collapse from any initial state.
\vskip 3pt
{\bf Theorem 22.} {\it Suppose
 $$l_1\ge 2d,\  l_1\le n-2d,\ n-2d<  l_2\le n-d.$$
 Then the system results in a state of collapse from any initial state.
 \vskip 3pt
 Proof.}  In accordance with Lemma 1, the system cannot be in a state of free movement.

 If, on a limit cycle, the system is not in neither a state of free motion nor collapse, then, in accordance with Lemma~3, on the limit cycle, the system results in at least one of the states  $(l_1+d,0),$ $(0,l_2+d),$ $(l_1,d),$ $(d,l_2).$
We shall prove that the system results in a state of collapse from any initial state.

Under the assumption that
  $$A(t_0)=(l_1+d,0),$$
 we have
  $$A(t_0+n-l_1-d)=(0,n-l_1-d),\ A(l_2+d)=(0,l_2+d),$$
  $$A(t_0+l_2+d)=(0,l_2+d),\ A(t_0+n)=(n-l_2-d,0),$$
  $$A(t_0+l_2+2d)=(d,l_2+2d-n),\ A(t_0+n+d)=(d,d).$$
Thus, from the states $(l_1+d_1, 0),$ $(0,l_2+d_2),$ the state results in a state of collapse.

Suppose
 $$A(t_0)=(l_1,d).$$
 Then we have
 $$A(t_0+n-l_1)=(0,n-l_1+d),\ A(t_0+l_2)=(0,l_2+d).$$
Therefore, from the state $(l_1,d),$ the system results in the state $(0,l_2+d),$ and, from this state, the system results in a state of collapse.

Assume that
 $$A(t_0)=(d,l_2).$$
 We have
 $$A(t_0+n-l_1)=(0,n-l_1+d),\ A(t_0+l_2)=(0,l_2+d).$$
Thus, from the state $(d,l_2),$ the system results in the state $(l_1+d,0),$ and, from this state, the system results in a state of collapse.
 Theorem~22 has been proved.

\vskip 3pt
 The condition of Theorem 22 holds, if, for example, $n=12,$ $l_1=5,$ $l_2=9,$ $d=2.$

 In this case,
 $$(7,0)\to (8,1)\to (9,2)\to (10,3)\to (11,4)\to (0,5)\to $$
 $$\to (0,6)\to (0,7)\to (0,8)\to (0,9)\to (0,10)\to $$
 $$\to (1,11)\to (2,0)\to (3,1)\to (3,2)\to (3,3),$$
 $$(6,2)\to (7,3)\to (8,4)\to (9,5)\to (10,6)\to (11,7)\to $$
 $$\to (0,8)\to (0,9),$$
 $$(2,7)\to (3,8)\to (4,9)\to (5,10)\to (6,11)\to (7,0)\to $$
 $$\to (8,0).$$

 \subsection*{4.23. Resulting in a state of collapse in the case
  $l_1\ge 2d,$ $n-2d< l_1 \le n-2d,$ $n-2d<  l_2\le n-d$
  }
\hskip 18pt If the condition of the following theorem holds, then the system results in a state of collapse from any initial state.
\vskip 3pt
{\bf Theorem 23.} {\it If the condition
 $$ n-2d< l_1  \le  l_2\le n-d,$$
 holds, then the system results in a state of collapse from any initial state.
 \vskip 3pt
 Proof.}  In accordance with Lemma 1 the system cannot be in a state of free movement.

 If, on a limit cycle, the system is not in neither a state of free motion nor collapse, then, in accordance with Lemma~3, on the limit cycle, the system results in at least one of the states  $(l_1+d,0),$ $(0,l_2+d),$ $(l_1,d),$ $(d,l_2).$
We shall prove that the system results in a state of collapse from any initial state.

 Assume that
  $$A(t_0)=(l_1+d,0).$$
 In this case,
  $$A(t_0+n-l_1-d)=(0,n-l_1-d),\ A(l_2+d)=(l_1+2d-n,d),\ A(t_0+n-l_1)=(d,d).$$
 Whence, from the state $(l_1+d, 0),$ the system results in a state of collapse.

 It is proved analogously that the system results in a state of collapse from the state $(0,l_2+d).$

 If
  $$A(t_0)=(l_1,d),$$
  then the following sequence of states is realized
  $$A(t_0+n-l_1-d)=(0,n-l_1-d),\ A(l_2+d)=(l_1+2d-n,d),\ A(t_0+n-l_1)=(d,d).$$
 Therefore, from the state $(l_1+d, 0),$ the system results in a state of collapse.

 Suppose
  $$A(t_0)=(l_1,d).$$
  Then we have
  $$A(t_0+n-l_1)=(0,n-l_1+d),\ A(t_0+l_2)=(0,l_2+d).$$
  Whence, from the state $(l_1,d),$ the system results in the state $(0,l_2+d),$ from which the system results in a state of collapse.

From the state $(d,l_2), $ the system results in the state
$(l_1+d,0),$ from which the system results in a state of collapse.
Theorem~23 has been proved.

 The condition of Theorem~23 holds, for example, under the assumptions that  $n=14,$ $l_1=9,$ $l_2=10,$ $d=3.$

 In this case, we have the following sequences of states
 $$(12,0)\to (13,1)\to (0,2)\to (1,3)\to (2,3)\to (3,3),$$
 $$(0,13)\to (1,0)\to (2,1)\to (3,2)\to (3,3),$$
 $$\to (1,11)\to (2,0)\to (3,1)\to (3,2)\to (3,3),$$
 $$(9,3)\to (10,4)\to (11,5)\to (0,6)\to (0,9)\to (0,10)\to $$
 $$\to (0,11)\to (0,12)\to (0,13),$$
 $$(3,10)\to (4,11)\to (5,12)\to (6,13)\to (7,0)\to (8,0)\to $$
 $$\to (9,0)\ to (10,0)\to (11,0)\to (12\to 0)\to (13,0).$$

 \subsection*{4.24. Resulting in a state of collapse in the case
 $n-2d< l_1 \le n-d,$ $l_2>n-d$
  }
\hskip 18pt If the condition of the following theorem holds, then the system results in a state of collapse from any initial state.
\vskip 3pt
{\bf Theorem 24.} {\it If the condition
 $$n-2d< l_1\le n-d,\  l_2> n-d,$$
 holds, then the system results in a state of collapse from any initial state.
 \vskip 3pt
 Proof.}  In accordance with Lemma 1 the system cannot be in a state of free movement.

 If, on a limit cycle, the system is not in neither a state of free motion nor collapse, then, in accordance with Lemma~3, on the limit cycle, the system results in at least one of the states  $(l_1+d,0),$ $(0,l_2+d),$ $(l_1,d),$ $(d,l_2).$
We shall prove that the system results in a state of collapse from any initial state.

 Suppose
  $$A(t_0)=(l_1+d,0).$$
  In this case,
  $$A(t_0+n-l_1-d)=(0,n-l_1-d),\ A(t_0+n-d)=(l_1-d,d),\ A(t_0+n-l_1)=(d,d).$$
 Whence, from the state $(l_1+d, 0),$ the system results in a state of collapse.

 Assume that
  $$A(t_0)=(0,l_2+d-n).$$
 We have
  $$A(t_0+n-l_2)=(n-l_2,d),\ A(t_0+d)=(d,d).$$
 Hence, from the state $(0,l_2+d-n),$ the system results in a state of collapse.

 Suppose
  $$A(t_0)=(l_1,d).$$
  Then we have
  $$A(t_0+n-l_1)=(0,n-l_1+d),\ A(l_2+d)=(l_1+2d-n,d),$$
  $$A(t_0+n-l_1)=(d,d).$$
 Hence, from the state $(l_1+d, 0),$ the system results in a statte of collapse.

 If
  $$A(t_0)=(0,l_2+d-n),$$
  then we have
  $$A(t_0+n-l_2)=(n-l_2,d),\ A(t_0+l_2)=(d,d).$$
Hence, from the state $(0,l_2+d-n),$ the state results in a state of collapse.

  Assume that
  $$A(t_0)=(l_1,d).$$
    Then we have
  $$A(t_0+n-l_1)=(0,n-l_1+d),\ A(t_0+n-l_1)=(0,n-l_1+d),$$
  $$A(t_0+l_2)=(0,l_2+d-n),$$
  Hence, from the state $(l_1,d),$  the system results in the state $(0,l_2+d-n),$ from which the system results in a state of collapse. Theorem~24 has been proved.

The condition of Theorem 24 holds if, for example, $n=12,$ $l_1=9,$ $l_2=11,$ $d=2.$

 The following sequences of states are realized
 $$(11,0)\to (0,1)\to (1,2)\to (2,2),$$
 $$(0,1)\to (1,2)\to (2,2),$$
 $$(1,11)\to (2,0)\to (3,1)\to (3,2)\to (3,3),$$
 $$(9,2)\to (10,3)\to (11,4)\to (0,5)\to (0,6)\to (0,7)\to $$
 $$\to (0,8)\to (0,9)\to (0,10)\to (0,11)\to (0,0)\to (0,1),$$
 $$(2,11)\to (3,0)\to (4,0)\to (5,0)\to (6,0)\to (7,0)\to $$
 $$\to (8,0)\to (9,0)\to (10,0)\to (11\to 0).$$

 \subsection*{4.25. Resulting in a state of collapse under the assumptions that
  $l_1>n-d,$ $l_2>n-d$
  }
\hskip 18pt If the condition of the following theorem holds, then the system results from a state of collapse from any initial sttate.
\vskip 3pt
{\bf Theoremc 25.} {\it If
 $$l_1>n-d,\  l_2>n-d,$$
 then the system results in a state of collapse from any iniitial state.
 \vskip 3pt
 Proof.} If the condition of Lemma 1 holds, then the system cannot be in a state of free motion..

 If, on a limit cycle, the system is not in either a state of free motion or collapse, then, in accordance of Lemma~3, on the limit cycle,  the system results in at least one of states $(l_1+d-n,0),$ $(0,l_2+d-n),$ $(l_1,d),$ $(d,l_2).$ We shall result in a state of collapse from any initial state.
 Assume that
  $$A(t_0)=(l_1+d-n,0).$$

 In this case,
  $$A(t_0+n-l_1)=(d,n-l_1),\ A(l_2+d)=(d,d).$$
 Whence, from the state $(l_1+d, 0),$ the system results in a state of collapse..

 It is proved analogously that the system results in the state of collapse from the state $(0,l_2+d-n).$

 Suppose
  $$A(t_0)=(l_1,d)$$
  In this case, we have
  $$A(t_0+n-l_1)=(0,n-l_1-d),\ A(l_2+l_2)=(0,l_2+d-n).$$
 Hence, from the state $(0,l_2+d-n),$ the system results in the state $(0,l_2+d-n),$ from which the system results in a state of collapse.

 Suppose
  $$A(t_0)=(l_1,d).$$
  Then we have
  $$A(t_0+n-l_1)=(0,n-l_1+d),\ A(t_0+l_2)=(0,l_2+d-n).$$
  Thus, from the state $(l_1,d),$ the state results in the state  $(0,l_2+d-n),$ from which the system results in a state of collapse.

 It is proved analogously that, from the state $(d,l_2),$ the state results in a state $1+d-n,0),$
 from which the system results in a state of collapse. Theorem~25 has been proved.

  The condition of Theorem 25 holds, for example, for $n=10,$ $l_1=8,$ $l_2=9,$ $d=3.$

 In this case, we have the following sequences of states
 $$(1,0)\to (2,1)\to (3,2)\to (3,3),$$
 $$(0,2)\to (1,3)\to (2,3)\to (3,3),$$
 $$(8,3)\to (9,4)\to (0,5)\to (0,6)\to (0,7)\to (0,8)\to $$
 $$\to (0,9)\to (0,0)\to (0,1)\to (0,2),$$
 $$(3,9)\to (4,0)\to (5,0)\to (6,0)\to (7,0)\to (8,0)\to $$
 $$\to (9,0)\ \to (0,0)\to (1,0).$$
 \vskip 3pt

 \section*{5. Conclusion}

 \hskip 18pt We consider different scenarios of behavior of two-contours system with different lengths of clusters. The number of these scenarios is sufficiently greater than in the case of clusters of the same length. In the later case, there 5 scenarios, [22]. In the former case, the scenarios are the following, Fig.~2--5.

\begin{figure}[ht!]
\centerline{\includegraphics[width=300pt]{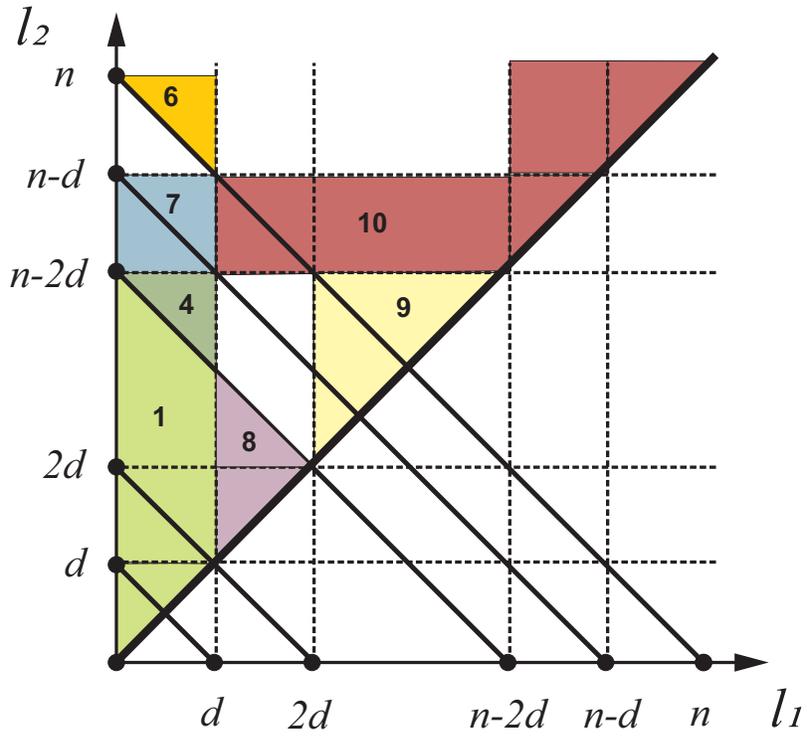}}
\caption{Scenarios of system behavior for different $0<l_1\le l_2<n$ and
$d<\frac{n}{4}$}
\end{figure}

\begin{figure}[ht!]
\centerline{\includegraphics[width=300pt]{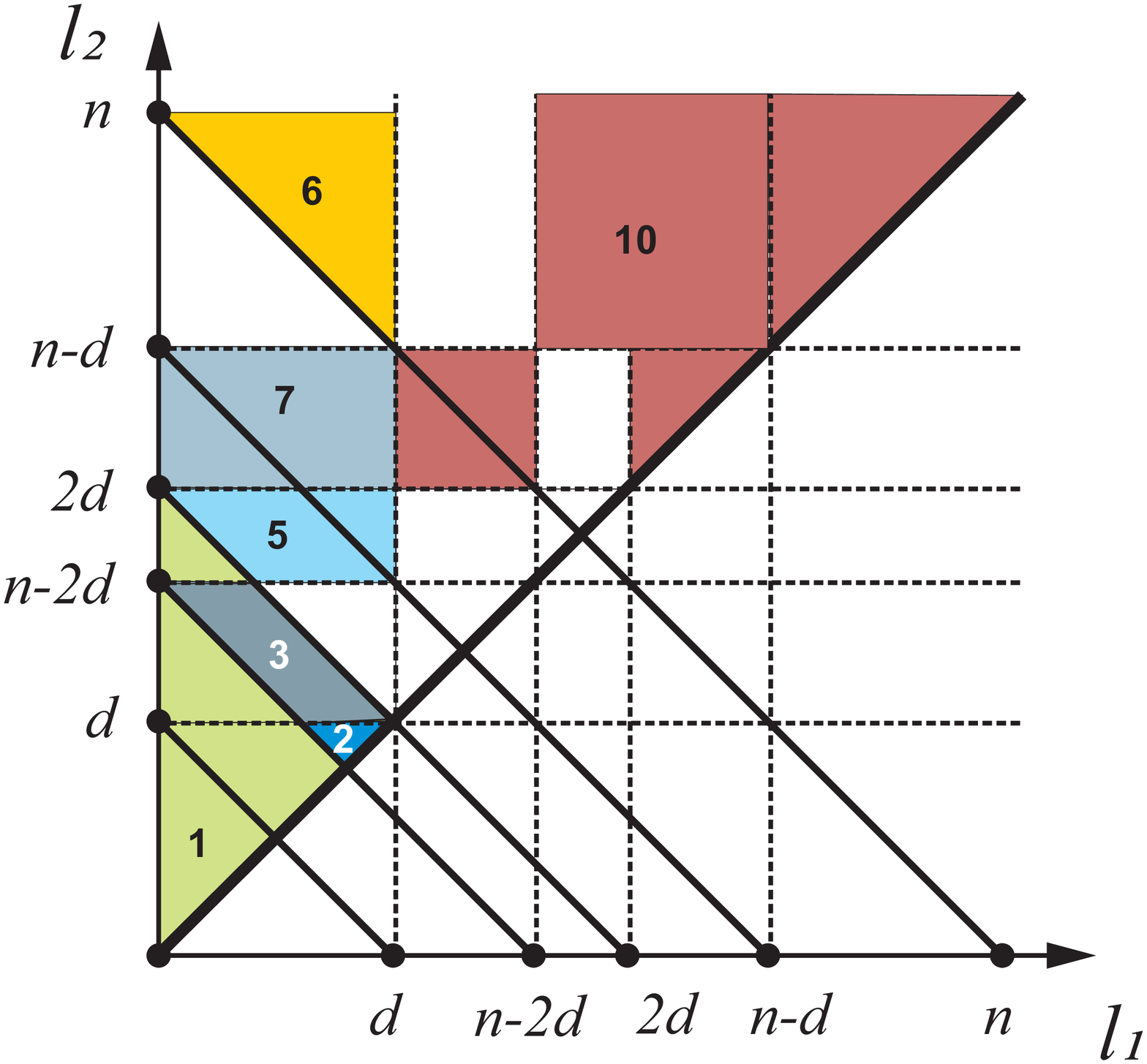}}
\caption{Scenarios of system behavior for different $0<l_1\le l_2<n$ and $\frac{n}{4}<d<\frac{n}{3}$}
\end{figure}

\begin{figure}[ht!]
\centerline{\includegraphics[width=300pt]{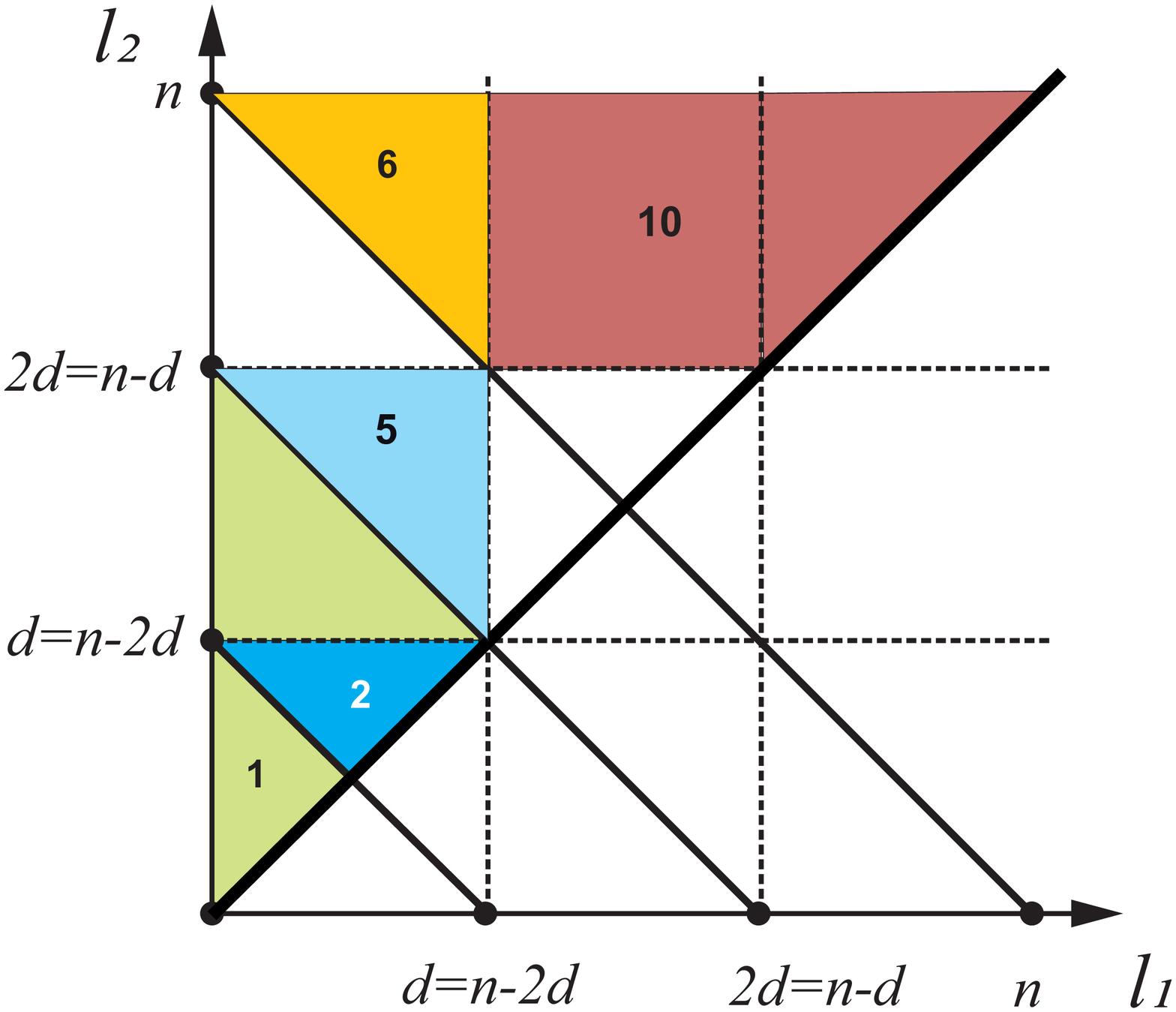}}
\caption{Scenarios of system behavior for different $0<l_1\le l_2<n$ and $d=\frac{n}{3}$}
\end{figure}

\begin{figure}[ht!]
\centerline{\includegraphics[width=300pt]{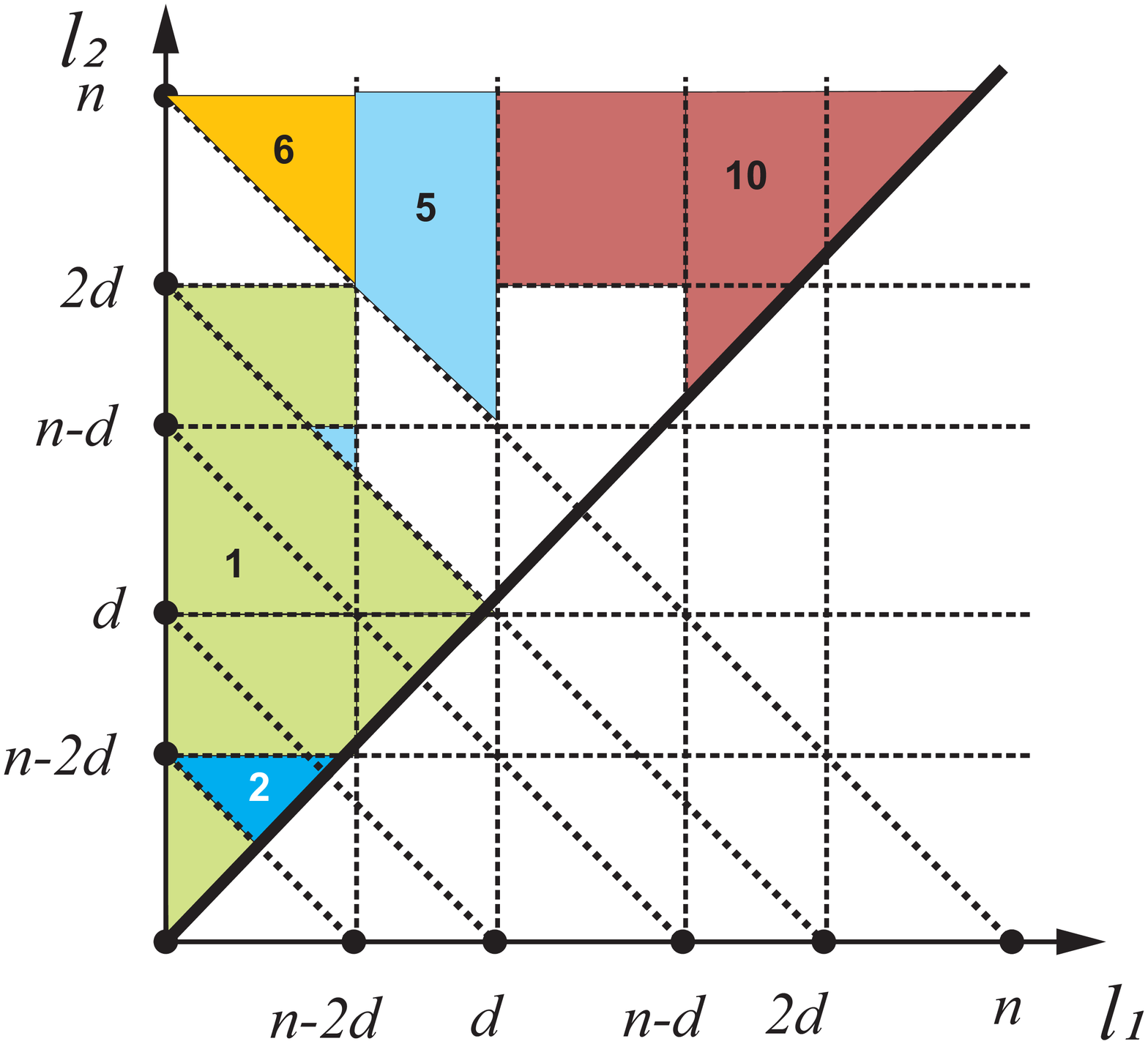}}
\caption{Scenarios of system behavior for different $0<l_1\le l_2<n$ and  $\frac{n}{3}<d<\frac{n}{2}$}
\end{figure}

\vskip 3pt
1) The system results in a state of free motion from any initial state.
\vskip 3pt
2) Depending on the initial state, $v=1$ or
$$v=\frac{n}{l_1+l_2+2d}.$$
\vskip 3pt
3) The average speed of each cluster is equal to
 $$v=\frac{n}{l_1+l_2+2d}$$
or
 $$v=\frac{n}{l_1+l_2+n-2d}.$$
\vskip 3pt
4) The average speed of clusters is equal to
 $$v=\frac{n}{l_1+l_2+2d}.$$
\vskip 3pt
5) The average speed of clusters is equal to
$$v=\frac{n}{l_1+l_2+n-2d}.$$
\vskip 3pt
6) The average speed of the cluster~1 is equal to
$$v_1=\frac{n}{2(l_1+l_2)},$$
and the average speed of the cluster~2 is equal to
 $$v_2=\frac{2}{l_1+l_2}.$$
\vskip 3pt
7) The average speed of the cluster~1 is equal to
$$v_1=\frac{1}{2},$$
and the cluster~2 is equal to
 $$v_2=1.$$
\vskip 3pt
8) Depending on the initial state, the system results in the state of free motion or collapse.
\vskip 3pt
9) The average speed equals 0 or
$$v=\frac{n}{l_1+l_2+2d}.$$
\vskip 3pt
10) The system results in a state of collapse from any initial state.

\section*{References}

1.  Nagel K.,  Schreckenberg M.A. Cellular automation models
for freeway traffic. {\it J. Phys. I.} 1992; 2(12):2221-2229.
DOI: 10.1051/jp1.1992277.\\
\vskip 5pt
2. Wolfram S., 1983, Statistical mechanics of cellular automata. {\it Rev. Mod. Phys.} 1983; 55: 601-644. DOI 10.1103/RevModPhys.55.601\\
\vskip 5pt
3. Belitzky V., Ferrary P.A. Invariant measures and
convergence properties for cellular automation 184 and related processes.
{\it J. Stat. Phys.} 2005; 118(3): 589-623.
DOI: 10.1007/s10955-044-8822-4.\\
\vskip 5pt
4. Blank M.L. Exact analysis of dynamical systems arising in
models of traffic flow. {\it Russian Mathematical Surveys.} 2000;
55:3, 562-563.\\
DOI: 10.4213/rm95
5. Gray L. and Grefeath D. The ergodic theory
of traffic jams. {\it J. Stat. Phys.} 2001; 105(3/4):413--452.
DOI: 10.1023/A:1012202706850.\\
\vskip 5pt
6. Kanai M, Nishinary K. and Tokihiro T. Exact solution and
asymptotic behavior of the asymmetric simple exclusion
process on a ring. {\it arXiv.0905.2795v1 [cond-mat-stat-mech] 18 May 2009.}\\
\vskip 5pt
7. Yashina M.V., Tatashev A.G. Traffic model based on synchronous and asynchronous exclusion processes.
Mathematical Methods in the Applied Sciences.  2020, vol.~43, issue 14, pp.~8136--8146.\\
Mathematical Methods in the Applied Sciences. First Published 04 February 2020.\\
DOI: 10.1102/mma6237
\vskip 5pt
8. 	Biham O., Middleton AA, Levine D. Self-organization and a dynamical transition in
traffic-flow models. {\it Phys. Rev. A.} 1992; 46(10):R6124-R6127.
DOI: 10.1003/PhysRevA.46.R6124.\\
\vskip 5pt
9. D'Souza R.M. Coexisting phases and lattice dependence of a cellular
automaton model for traffic flow. {\it Phys. Rev. E.} 2005;
71(6):066112. DOI: 10.1103/PhysRevE.71.066112.\\
\vskip 5pt
10. Angel O., Holroyd AE, Martin JB.
The Jammed Phase of the Biham–Middleton–Levine Traffic Model.
{\it Electronic Communications in Probability.} 2005; 10:167–178.\\
DOI: 10.1214/ECP.v10-1148.
\vskip 5pt
11. Austin T., Benjamini I. {\it For what number of cars must self-organization occur
in the Biham-Middleton-Levine traffic model from any possible starting configuration.}
2006; arXiv.math/0607759.\\
\vskip 5pt
12. Bugaev AS, Buslaev AP,
Kozlov V.V, Yashina M.V. Distributed problems of monitoring
and modern approaches to traffic modeling,
{\it 14th International IEEE Conference
on Intelligent Transactions Systems (ITSC~2011),
Washington, USA, 5--7.10.2011.} 2011; 477--481.\\
\vskip 5pt
13.  Kozlov VV, Buslaev AP, Tatashev A.G.
On synergy of totally connected flow on
chainmails. {\it CMMSE-2013, Cadis, Spain}, 2013; 3, 861--873.
\vskip 5pt
14.  Buslaev A.P., Fomina M.Yu., Tatashev A.G., Yashina M.V.
On discrete flow networks model spectra:
statement, simulation, hypotheses.
{\it J. Phys.: Conf. Ser.} 2018. 1053(012034).
DOI: 10.1088/1742/6596/1053/1/012034.\\
\vskip 5pt
15.  Buslaev A.P.,  Tatashev A.G. Flows on discrete
traffic flower. {\it Journal of Mathematics Research}. 2017.
9(1), 98--108. DOI 10.5539/jmr.v9n1p98.\\
\vskip 5pt
16. Buslaev AP, Tatashev AG. Exact results for discrete
dynamical systems on a pair of contours. {\it Math.
Meth. Appl. Sci.} 2018; 41(17):7283-7294.\\
\vskip 5pt
19. Buslaev A.P., Tatashev A.G, Yashina M.V. Flows spectrum on closed trio of contours
{\it Eur. J. Pure Appl. Math,} 2018.
{\bf 11(3),} 893-897.\\ (DOI 10.29020/nybg.ejpam.v11i1.3201)
\vskip 5pt
21. Yashina M.V., Tatashev A.G. Spectral cycles and average velocity of clusters in discrete two-contours system with two nodes. Math. Method. Appl. Sci., 2020, vol.~43, issue7, pp.~4303--4316. DOI: 10.1102/mma6194
\vskip 3pt
22. Yashina M., Tatashev A. Uniform cluster traffic model on closed two-contours system with two non-symmetrical common nodes. Traffic and Granular Flow '19. (In print.)
\vskip 3pt
23.  Yashina M.V., Tatashev A.G.,  Fomina M.Y. Optimization of velocity mode in Buslaev two-contour networks via competition resolution rules. International Journal of Interactive Mobile Technologies, 2020, vol.~14, no.~10, pp.~61--73.\\
DOI: 103991/ijim.v14i10.14641

\end{document}